\documentclass[10pt]{amsart}

\input{diagrams}

\newtheorem{proposition}{Proposition}[section]
\newtheorem{lemma}[proposition]{Lemma}
\newtheorem{corollary}[proposition]{Corollary}
\newtheorem{theorem}[proposition]{Theorem}

\theoremstyle{definition}
\newtheorem{definition}[proposition]{Definition}

\theoremstyle{remark}
\newtheorem{remark}[proposition]{Remark}
\newtheorem{remarks}[proposition]{Remarks}

\newcommand{\thlabel}[1]{\label{th:#1}}
\newcommand{\thref}[1]{Theorem~\ref{th:#1}}
\newcommand{\selabel}[1]{\label{se:#1}}
\newcommand{\seref}[1]{Section~\ref{se:#1}}
\newcommand{\lelabel}[1]{\label{le:#1}}
\newcommand{\leref}[1]{Lemma~\ref{le:#1}}
\newcommand{\prlabel}[1]{\label{pr:#1}}
\newcommand{\prref}[1]{Proposition~\ref{pr:#1}}
\newcommand{\colabel}[1]{\label{co:#1}}
\newcommand{\coref}[1]{Corollary~\ref{co:#1}}
\newcommand{\relabel}[1]{\label{re:#1}}
\newcommand{\reref}[1]{Remark~\ref{re:#1}}

\newcommand{\delabel}[1]{\label{de:#1}}
\newcommand{\deref}[1]{Definition~\ref{de:#1}}
\newcommand{\eqlabel}[1]{\label{eq:#1}}
\newcommand{\equref}[1]{(\ref{eq:#1})}

\newcommand{\Hom}{{\rm Hom}}
\newcommand{\End}{{\rm End}}

\newcommand{\Ker}{{\rm Ker}\,}

\newcommand{\im}{{\rm Im}\,}

\def\Ab{\underline{\underline{\rm Ab}}}

\def\ot{\otimes}

\def\CC{{\mathbb C}}
\def\ZZ{{\mathbb Z}}
\def\QQ{{\mathbb Q}}

\def\DD{{\mathbb D}}
\def\TT{{\mathbb T}}
\def\SS{{\mathbb S}}

\newcommand{\Cc}{\mathcal{C}}
\newcommand{\Dd}{\mathcal{D}}

\newcommand{\Mm}{\mathcal{M}}

\def\*C{{}^*\hspace*{-1pt}{\Cc}}
\def\text#1{{\rm {\rm #1}}}

\def\rightrightarrows{{\pile{\longrightarrow\\\longrightarrow}}}

\def\ol{\overline}
\def\ul{\underline}
\def\dul#1{\underline{\underline{#1}}}

\begin{document}
\title[Galois theory  for comatrix corings ]{Galois theory for comatrix corings:
descent theory, Morita theory, Frobenius and separability properties}
\author{S. Caenepeel}
\address{Faculty of Applied Sciences,
Vrije Universiteit Brussel, VUB, B-1050 Brussels, Belgium}
\email{scaenepe@vub.ac.be}
\urladdr{http://homepages.vub.ac.be/\~{}scaenepe/}
\author{E. De Groot}
\address{Faculty of Applied Sciences,
Vrije Universiteit Brussel, VUB, B-1050 Brussels, Belgium}
\email{edegroot@vub.ac.be}
\urladdr{http://homepages.vub.ac.be/\~{}edegroot/}
\author{J. Vercruysse}
\address{Faculty of Applied Sciences,
Vrije Universiteit Brussel, VUB, B-1050 Brussels, Belgium}
\email{joost.vercruysse@vub.ac.be}
\urladdr{http://homepages.vub.ac.be/\~{}jvercruy/}
\thanks{}


\begin{abstract}
El Kaoutit and G\'omez Torrecillas introduced comatrix corings, generalizing
Sweedler's canonical coring, and proved a new version of the Faithfully Flat Descent
Theorem. They also introduced Galois corings, as corings isomorphic to a
comatrix coring. In this paper, we further investigate this theory. We prove a
new version of the Joyal-Tierney Descent Theorem, and generalize the Galois Coring
Structure Theorem. We associate a Morita context to a coring with a fixed comodule,
and relate it to Galois-type properties of the coring. An affineness criterion is
proved in the situation where the coring is coseparable. Further properties of the
Morita context are studied in the situation where the coring is (co)Frobenius.
\end{abstract}

\maketitle

\section*{Introduction}
Corings were introduced by Sweedler in \cite{Sweedler65}. Takeuchi \cite{Tak}
remarked that entwined modules introduced in \cite{BrzezinskiM}
can be viewed as examples of comodules over a coring. Takeuchi's observation
has caused a revival of the theory of corings: it became clear that a number of
results from Hopf algebra and related areas can be at the same time reformulated
and generalized using the language of corings. The computations become
simpler, more natural and more transparent. Graded modules, Hopf modules,
Long dimodules, Yetter-Drinfeld modules, entwined modules and weak entwined modules
are special cases of comodules over a coring. Corings can be used to study
properties of functors between categories of graded modules, Hopf modules,...
This was discussed by Brzezi\'nski in \cite{Brzezinski02},
the first of a series of papers illustrating the importance of corings. For a complete
list of references, we refer to the recent monograph \cite{BrzezinskiWisbauer},
in which a number of applications of corings are presented.

Corings can be used to present an elegant presentation of descent and Galois theory.
The idea appears already in \cite{Brzezinski02}, and was further investigated in
\cite{Abu,Caenepeel03,CVW,W}. Given a ring morphism $B\to A$, one can introduce
the category of descent data, see for example \cite{KO} in the case where $A$ and $B$
are commutative, and \cite{Cipolla} in the noncommutative case. A descent datum
turns out to be a comodule over the Sweedler canonical coring $\Dd=A\ot_BA$.
A Galois coring is then by definition a coring that is isomorphic to the canonical coring,
and a Galois descent datum is a comodule over this coring. For example, if
$H$ is Hopf algebra, and $A/B$ is an $H$-Galois extension in the sense of
\cite{Schneider90}, then $A\ot H$ can be made into a coring over $A$, which is
isomorphic to the canonical coring $A\ot_BA$. In a similar way, classical Galois
extensions, strongly graded rings, and coalgebra Galois extensions (see
\cite{BrzezinskiH}) can be introduced using Galois corings.

In \cite{Kaoutit}, El Kaoutit and G\'omez Torrecillas look at a more general version
of the descent problem: in the classical situation, we take  a ring morphism $B\to A$,
and try to descend modules defined over $A$ to modules defined over $B$.
This theory can be generalized to the situation where $A$ and $B$ are connected
by a $(B,A)$-bimodule $\Sigma$. The descent data are now comodules over
the {\it comatrix coring}, which is equal to $\Sigma^*\ot_B\Sigma$ as an
$A$-bimodule. El Kaoutit and G\'omez Torrecillas prove the faithfully flat descent
theorem in this setting, and introduce a generalized notion of Galois coring;
basically it is a coring that is isomorphic to a comatrix coring. Comatrix corings have
been studied also in \cite{Br3,BrzezinskiG}.

In this paper, we further investigate this theory. In \seref{2}, we look at descent theory.
The most famous, but not most general result in the classical setting is
the faithfully flat descent theorem: if $A/B$ is faithfully flat, then the category of
descent data (comodules over the canonical coring) is equivalent to the category
of $B$-modules. A more general result, due to Joyal and Tierney (unpublished)
is the following: if $A$ and $B$ are commutative, then we have the desired equivalence
if and only if $i:\ B\to A$ is pure as a map of $B$-modules. We will present a 
generalization of the Joyal-Tierney Theorem to the comatrix coring situation:
a sufficient condition for category equivalence is now that $i:\ B\to \End_A(\Sigma)$
is pure as map of left $B$-modules, and that $i$ maps $B$ into the center of $\End_A(\Sigma)$.
Our proof is inspired by Mesablishvili's proof of the Joyal-Tierney Theorem.

In \seref{3}, we recall the definition of Galois coring from \cite{Kaoutit}; we can directly
translate some of the results of \seref{2}. The main results of
the Section are Theorems \ref{th:3.8} and \ref{th:3.9}, which are generalizations of
the Galois Coring Structure Theorem from \cite{W}.

In \seref{4}, we associate a Morita context to a comodule over a coring. It can be
viewed as a dual version of the classical Morita context associated to a module
over a ring. Actually, there is morphism from our Morita context to the Morita context
associated to $\Sigma$ viewed as a module over the dual coring, and these are
isomorphic under some finiteness assumptions. We can apply the Morita context
to obtain more equivalent conditions for the Galois descent in the situation where
the coring is finitely generated and projective as an $A$-module (see \thref{4.11}).

A coring is Galois if a certain map (called the canonical map) from the canonical
coring to the coring is bijective. Sometimes surjectivity is sufficient; classical results
in the Hopf algebra case are in \cite{Schneider90}. These results were improved
recently in \cite{SchSch}; in the case of Doi-Hopf modules, some results were
presented in  \cite{Militaru}. In \seref{5}, we give a result of this type in the coring
situation: surjectivity is sufficient in the situation where $\Cc$ is a coseparable
coring.

The Morita context that we introduce in \seref{4} is in fact a generalization
of a Morita context introduced by Doi \cite{D}. Morita contexts similar to the
one of Doi were studied by Cohen, Fischman and Montgomery in \cite{CF} and \cite{CFM}.
These are different from the one of Doi, in the sense that the two connecting modules
in the context are equal to the underlying algebra $A$. On the other hand, they are
more restrictive, in the sense that they only work for finite dimensional Hopf algebras
over a field (see \cite{CF}) or Frobenius Hopf algebras over a commutative ring
(see \cite{CFM}). This has been clarified in \cite{CVW}, using the notion of Frobenius
coring. In \seref{6}, we study the Morita context associated to a Frobenius coring
with a fixed comodule $\Sigma$. It turns out that the connecting modules in the context
are then precisely $\Sigma$ and its right dual $\Sigma^*$; in the case where
$\Sigma=A$, the situation studied in \cite{CVW}, the two connecting modules are
then isomorphic to $A$. Weaker results are obtained in the situation where
$\Cc$ is coFrobenius.

It is well-known that the set of right $\Cc$-comodule structures on $A$ corresponds bijectively to the set of grouplike elements of the coring $\Cc$. As we already
indicate, if we take $\Sigma=A$, then we recover the ``classical" Galois theory for corings. Another possible choice is $\Sigma=\Cc$, at least in the case where
$\Cc$ is finitely generated and projective as a right $A$-module. This situation
is examined in \seref{7}.

\section{Preliminary results}\selabel{1}
Let $A$ be a ring. Recall that an $A$-coring is a comonoid in the monoidal
category ${}_A\Mm_A$. Thus a coring is a triple $(\Cc,\Delta_{\Cc},\varepsilon_{\Cc})$,
where $\Cc$ is an $A$-bimodule, and $\Delta:\ \Cc\to \Cc\ot_A\Cc$ and
$\varepsilon_\Cc:\ \Cc\to A$ are $A$-bimodule maps such that
$$(\Delta_{\mathcal C}\ot_A  {\mathcal C}) \circ \Delta_{\mathcal
C}=({\mathcal C}\ot_A \Delta_{\mathcal C}) \circ
\Delta_{\mathcal C}, \quad (\varepsilon_{\mathcal C} \ot_A
{\mathcal C}) \circ \Delta_{\mathcal C}= (Id_{\mathcal C} \ot_A
\varepsilon_{\mathcal C}) \circ \Delta_{\mathcal C}={\mathcal C}.$$
$\Delta_\Cc$ is called the comultiplication, and $\varepsilon_\Cc$ is called
the counit. We use the Sweedler-Heyneman notation
$$\Delta_{\Cc}= c_{(1)}\ot_A c_{(2)},$$
where the summation is implicitely understood. A right $\Cc$-comodule
is a couple $(M,\rho^r)$, where $M$ is a right $A$-module, and
$\rho^r:\ M\to M\ot_A\Cc$ is a right $A$-linear map, called the coaction,
satisfying the conditions
$$(\rho^r \ot_A Id_{\mathcal C}) \circ \rho^r=(Id_M \ot_A
\Delta_{\mathcal C}) \circ \rho^r, \quad (Id_M \ot_A
\varepsilon_{\mathcal C}) \circ \rho^M = Id_M.$$
We use the following Sweedler-Heyneman notation for right coactions:
$$\rho^r(m)= m_{[0]}\ot_Am_{[1]}.$$
Let $M$ and $N$ be right $\Cc$-comodules. A right $R$-linear map
$f:\ M\to N$ is called right $\Cc$-colinear if
$$\rho^N(f(m))=f(m_{[0]})\ot_A m_{[1]},$$
for all $m\in M$. The category of right $\Cc$-comodules and right $\Cc$-colinear
maps is denoted by $\Mm^\Cc$. The full subcategory consisting of right
$\Cc$-comodules that are finitely generated and projective as a right
$A$-module is denoted by $\Mm_{\rm fgp}^\Cc$.\\ 
In a similar way, we define left $\Cc$-comodules $(M,\rho^l)$, with
$\rho^l:\ M\to \Cc\ot_A M$ a left $A$-module map. The Sweedler-Heyneman notation
for left coactions is
$$\rho^l(m)= m_{[-1]}\ot_Am_{[0]}.$$
The category of left $\Cc$-comodules and left $\Cc$-colinear maps is denoted by
${}_{\rm fgp}^\Cc\Mm$.
Let $\Sigma\in \Mm_A$. Then $\Sigma^*\in {}_A\Mm$, with left $A$-action
$(af)(u)=af(u)$, for all $a\in A$ and $u\in \Sigma$. $\Sigma$ is finitely generated
and projective in $\Mm_A$ if and only if there exists a (unique) $e=\sum_i e_i\ot_Af_i\in
\Sigma\ot_A\Sigma^*$, by abuse of language called the dual basis of $\Sigma$, such that
\begin{equation}\eqlabel{1.1}
u=\sum_i e_if_i(u)~~{\rm and}~~f=\sum_i f(e_i)f_i,
\end{equation}
for all $u\in \Sigma$ and $f\in \Sigma^*$. 
In this case, $\Sigma^*$ is finitely generated projective in ${}_A\Mm$. We obtain
a pair of inverse  equivalences $((\bullet)^*,{}^*(\bullet))$
between $\Mm_{A,{\rm fgp}}$ and ${}_{A,{\rm fgp}}\Mm^{\rm op}$.

\begin{proposition}\prlabel{1.1}
Let $\Cc$ be an $A$-coring. We have a  pair of inverse equivalences between the
categories $\Mm_{\rm fgp}^\Cc$ and ${}_{\rm fgp}^\Cc\Mm$.
\end{proposition}

\begin{proof}
Take $(\Sigma,\rho^r)\in \Mm_{\rm fgp}^\Cc$, and let $e$ be a finite dual basis.
Consider
\begin{equation}\eqlabel{1.1.1}
\rho^l:\ \Sigma^*\to \Cc\ot_A\Sigma^*,~~\rho^l(f)=\sum_i f(e_{i[0]})e_{i[1]}\ot_A f_i.
\end{equation}
Let us show that $(\Sigma^*,\rho^l)\in {}_{A,{\rm fgp}}\Mm^{\rm op}$.
From \equref{1.1}, it follows that
$$u_{[0]}\ot_A u_{[1]}=\sum_i e_{i[0]}\ot_A e_{i[1]}f_i(u),$$
hence
$$f(u_{[0]})u_{[1]}=\sum_i f(e_{i[0]})e_{i[1]}f_i(u).$$
Using this property, we find
\begin{eqnarray*}
&&\hspace*{-2cm}
(I\ot_A\rho^l)(\rho^l(f))=
\sum_{i,j}f(e_{i[0]})e_{i[1]}\ot_A f_i(e_{j[0]})e_{j[1]}\ot_Af_j\\
&=&\sum_{i,j}f(e_{i[0]})e_{i[1]}f_i(e_{j[0]})\ot_A e_{j[1]}\ot_Af_j\\
&=&\sum_j f(e_{j[0]})e_{j[1]}\ot_A e_{j[2]}\ot_A f_j=
(\Delta\ot_A I)(\rho^l(f))
\end{eqnarray*}
and
$$(\varepsilon_\Cc\ot_A I)(\rho^l(f))=\sum_i f(e_{i[0]})\varepsilon_\Cc(e_{i[1]})f_i=f,$$
as needed. All the other verifications are straightforward and left to the reader.
\end{proof}

Now we consider a second ring $B$. We call $M$ a $(B,\Cc)$-bicomodule if
$M$ is a $(B,A)$-bimodule and a right $\Cc$-comodule such that
\begin{equation}\eqlabel{1.2.1}
\rho^r(bm)=bm_{[0]}\ot_A m_{[1]},
\end{equation}
for all $b\in B$ and $m\in M$. This means that the canonical map
$$l:\ B\to \End_A(M),~~l_b(m)=bm$$
factorizes through $\End^\Cc(M)$. The category of $(B,\Cc)$-bicomodules
and left $B$-linear right $\Cc$-colinear maps is denoted ${}_B\Mm^\Cc$.
The full subcategory consisting of $(B,\Cc)$-bicomodules that are finitely
generated and projective as right $A$-modules is denoted by
${}_B\Mm^\Cc_{\rm fgp}$. We will use a similar notation for left $\Cc$-comodules.\\

Let $\Cc$ be an $A$-coring, and consider $M\in \Mm^\Cc$ and $N\in {}^\Cc\Mm$.
$$M\ot^\Cc N=\{\sum_j m_j\ot n_j\in M\ot_A N~|~
\sum_j\rho^r( m_j)\ot n_j=\sum_j m_j\ot \rho^l(n_j)\}$$
is called the cotensor product of $M$ and $N$. Observe that it is the
equalizer of $\rho^r\ot_A I_N$ and $I_M\ot_A \rho^l$. $M\ot^\Cc N$ is
an abelian group, but in some cases it has more structure. The proof of the
following result is trivial.

\begin{lemma}\lelabel{1.2}
If $M\in {}_B\Mm^\Cc$ and $N\in {}^\Cc\Mm_D$, then $M\ot^\Cc N\in
{}_B\Mm_D$.
\end{lemma}

Also observe that
$M\ot_A\varepsilon_\Cc:\ M\ot^\Cc \Cc\to M$
is an isomorphism with inverse $\rho^r$. Another property in the same style is
the following:

\begin{lemma}\lelabel{1.2b}
Let $L\in \Mm^\Cc$, $M\in \Mm_A$. Then we have an isomorphism
$$\alpha:\ \Hom_A(L,M)\to \Hom^\Cc(L,M\ot_A\Cc),$$
given by
$$\alpha(f)(l)=f(l_{[0]})\ot_A l_{[1]}~~{\rm and}~~
\alpha^{-1}(\varphi)=(I_M\ot_A\varepsilon_\Cc)\circ\varphi.$$
\end{lemma}

\begin{proof}
It is clear that $\alpha(f)$ is right $\Cc$-colinear, for any $f\in \Hom_A(L,M)$.
Furthermore
$$(\alpha^{-1}(\alpha(f))(l)=
(M\ot_A\varepsilon_\Cc)(f(l_{[0]})\ot_A l_{[1]})=f(l).$$
Take $\varphi\in  \Hom^\Cc(L,M\ot_A\Cc)$. If $\varphi(l)=\sum_j m_j\ot_A c_j$,
then $\varphi(l_{[0]})\ot_A l_{[1]}=\sum_j m_j\ot_A c_{j[1]}\ot_A c_{j[2]}$, and
\begin{eqnarray*}
&&\hspace*{-2cm}
(\alpha\circ \alpha^{-1})(\varphi)(l)=
\alpha^{-1}(\varphi)(l_{[0]})\ot_A l_{[1]}\\
&=& (M\ot_A\varepsilon_\Cc)(\varphi(l_{[0]}))\ot_A l_{[1]}\\
&=& \sum_j m_j\ot_A \varepsilon_\Cc(c_{j[1]})\ot_A c_{j[2]}=\varphi(l).
\end{eqnarray*}
\end{proof}

Let $A$ and $B$ be rings, and $\Sigma\in {}_B\Mm^\Cc$. Then $\Sigma^*\in {}_A\Mm_B$ via
$$(afb)(u)=af(bu).$$
If $e$ is the dual basis, then $e\in (\Sigma\ot_A \Sigma^*)^B$. Indeed,
for all $b\in B$, we have
\begin{eqnarray*}
&&\hspace*{-2cm}
be=\sum_i be_i\ot_A f_i=\sum_{i,j} e_jf_j(be_i)\ot_A f_i\\
&=& \sum_{i,j} e_j\ot_Af_j(be_i) f_i= \sum_{i,j} e_j\ot_A(f_jb)(e_i) f_i\\
&=& \sum_j e_j\ot_A f_jb=eb
\end{eqnarray*}
We have a ring isomorphism
\begin{equation}\eqlabel{1.2.2}
(\bullet)^*:\ \End_A(\Sigma)\to {}_A\End(\Sigma^*)^{\rm op},
\end{equation}
 sending $f$ to its dual map $f^*$. It restricts to an isomorphism
 $$(\bullet)^*:\ \End^\Cc(\Sigma)\to {}^\Cc\End(\Sigma^*)^{\rm op},$$
 and we have
 $$r=(\bullet)^*\circ l:\ B\to {}^\Cc\End(\Sigma^*)^{\rm op},~~r_b(f)=fb.$$
 
\begin{proposition}\prlabel{1.3}
Let $\Cc$ be an $A$-coring, $M\in \Mm^\Cc$ and $\Sigma\in \Mm^\Cc_{\rm fgp}$.
Then the canonical isomorphism $\alpha:\ \Hom_A(\Sigma, M)\to
M\ot_A\Sigma^*$ restricts to an isomorphism
$$\Hom^\Cc(\Sigma, M)\cong M\ot^{\Cc}\Sigma^*.$$
\end{proposition}

\begin{proof}
Recall that $\alpha(\varphi)= \sum_i \varphi(e_i)\ot_A f_i$,
and $\alpha^{-1}(m\ot_A g)=\varphi$ with $\varphi(u)=mg(u)$,
for all $\varphi\in \Hom_A(\Sigma, M)$, $m\in M$, $g\in \Sigma^*$ and
$u\in \Sigma$.\\
Take $\sum_j m_j\ot_A g_j\in M\ot_A\Sigma^*$, and let
$\varphi=\alpha^{-1}(\sum_j m_j\ot_A g_j)\in \Hom_A(\Sigma,M)$ be the corresponding map.
Then $\varphi\in \Hom^\Cc(\Sigma,M)$ if and only if
\begin{equation}\eqlabel{1.3.2}
\varphi(u_{[0]})\ot_A u_{[1]}=\varphi(u)_{[0]}\ot_A\varphi(u)_{[1]},
\end{equation}
for all $u\in \Sigma$.
The right hand side amounts to
$$\varphi(u)_{[0]}\ot_A\varphi(u)_{[1]}=\sum_j m_{j[0]}\ot_A m_{j[1]}g_j(u)$$
and we compute the left hand side:
\begin{eqnarray*}
&&\hspace*{-2cm}
\varphi(u)_{[0]}\ot_A\varphi(u)_{[1]}=
\sum_j m_jg_j(u_{[0]})\ot_A u_{[1]}=
 \sum_{i,j} m_jg_j(e_{i[0]})\ot_A e_{i[1]}f_i(u)\\
&=&\sum_{i,j} m_j\ot_Ag_j(e_{i[0]}) e_{i[1]}f_i(u)=
\sum_j m_j\ot_A g_{j[-1]}g_{j[0]}(u).
\end{eqnarray*}
If \equref{1.3.2} holds, then we find for all $i$:
$$\sum_j m_j\ot_A g_{j[-1]}g_{j[0]}(e_i)
=\sum_j m_{j[0]}\ot_A m_{j[1]}g_j(e_i),$$
and consequently,
$$\sum_{i,j} m_j\ot_A g_{j[-1]}g_{j[0]}(e_i)\ot_A f_i
=\sum_{i,j} m_{j[0]}\ot_A m_{j[1]}g_j(e_i))\ot_A f_i,$$
or
$$\sum_{i,j} m_j\ot_A g_{j[-1]}\ot_A g_{j[0]}(e_i) f_i
=\sum_{i,j} m_{j[0]}\ot_A m_{j[1]}\ot_A g_j(e_i))f_i,$$
and, finally,
\begin{equation}\eqlabel{1.3.3}
\sum_j m_j\ot_A \rho^l(g_j)=\sum_j \rho^r(m_j)\ot_A g_j.
\end{equation}
Conversely, if \equref{1.3.3} holds, then \equref{1.3.2} follows after applying
the last tensor factor to $u\in \Sigma$.
\end{proof}

\begin{proposition}\prlabel{1.4}
Let $A$ and $B$ be rings, $\Cc$ an $A$-coring, and $\Sigma\in {}_B\Mm^\Cc_{\rm fgp}$.
Then we have the following two pairs of adjoint functors $(F,G)$ and $(F',G')$:
$$F:\ \Mm_B\to \Mm^\Cc,~~F(N)=N\ot_B\Sigma$$
$$G:\ \Mm^\Cc\to \Mm_B,~~G(M)=\Hom^\Cc(\Sigma,M)\cong M\ot^\Cc\Sigma^*$$
and
$$F':\ {}_B\Mm\to {}^\Cc\Mm,~~F'(N)=\Sigma^*\ot_BN$$
$$G':\ {}^\Cc\Mm\to {}_B\Mm,~~G'(M)={}^\Cc\Hom(\Sigma^*,M)\cong \Sigma\ot^\Cc M$$
\end{proposition}

\begin{proof}
We will only give the unit and counit of the first adjunction, leaving all other verifications
to the reader. For $N\in \Mm_B$:
$$\nu_N:\ N\to \Hom^\Cc(\Sigma,N\ot_B\Sigma),~~\nu_N(n)(u)=n\ot_B u,$$
or
$$\nu_N:\ N\to (N\ot_B\Sigma)\ot^\Cc\Sigma^*,~~
\nu_N(n)=\sum_i (n\ot_B e_i)\ot_A f_i,$$
and for $M\in\Mm^\Cc$:
$$\zeta_M:\ \Hom^\Cc(\Sigma,M)\ot_B \Sigma\to M,~~\zeta_M(\varphi\ot_Bu)=\varphi(u),$$
or
$$\zeta_M:\ (M\ot^\Cc\Sigma^*)\ot_B \Sigma\to M,~~\zeta_M((\sum_j m_j\ot_A g_j)\ot_Bu)=\sum_j m_jg_j(u).$$
\end{proof}

Our aim is to determine when $(F,G)$ and $(F',G')$ are inverse equivalences. 
We will first do this in the case where $\Cc$ is the so-called comatrix coring associated
to a bimodule.

\section{Comatrix corings and descent theory}\selabel{2}
Let $A$ and $B$ be rings, and $\Sigma\in {}_B\Mm_A$ a bimodule that is
finitely generated and projective as a right $A$-module, with finite
dual basis $e=\sum_i e_i\ot_A f_i$. 
Then $\Dd= \Sigma^*\ot_B \Sigma$ is an $A$-coring; comultiplication and
counit are given by the formulas
$$\Delta_\Dd:\  \Dd\to \Dd\ot_A\Dd,~~\Delta_\Dd(f\ot_B u)=f\ot_B e\ot_Bu;$$
$$\varepsilon_\Dd:\  \Dd\to A,~~\varepsilon_\Dd(f\ot_B u)=f(u).$$
$\Dd$ is called the comatrix coring associated to the bimodule $\Sigma$;
comatrix corings have been studied in \cite{Kaoutit} and \cite{BrzezinskiG}. 
Also $\Sigma$ is a right $\Dd$-comodule and $\Sigma^*$ is a left
$\Dd$-comodule; the coactions are given by the formulas
$$\rho^r(u)=e\ot_B u~~{\rm and}~~\rho^l(f)=f\ot_B e.$$
$\Sigma\in {}_B\Mm^\Dd$,
since \equref{1.2.1} holds: 
\begin{eqnarray*}
&&\hspace*{-2cm}
\rho^r(bu)=\sum_i e_i\ot_A f_i\ot_B bu=
\sum_i e_i\ot_A f_ib\ot_B u\\
&=& \sum_i be_i\ot_A f_i\ot_B u= bu_{[0]}\ot_A u_{[1]},
\end{eqnarray*}
for all $b\in B$ and $u\in \Sigma$, where we used the fact that $e\in (\Sigma\ot_A
\Sigma^*)^B$. For any $M\in \Mm^\Dd$, we have that
$$\Hom^{\Dd}(\Sigma, M)\cong M\ot^\Dd \Sigma^*,$$
the subspace of $\sum_j m_j\ot_A g_j\in M\ot_A\Sigma^*$ satisfying
$$\sum_j \rho^M(m_j)\ot_A g_j= \sum_j m_j\ot_A g_j\ot_B e.$$
In particular, 
$$T=\End^{\Dd}(\Sigma)\cong\{x\in \Sigma\ot_A\Sigma^*~|~e\ot_Bx=x\ot_B e\}.$$
Following \prref{1.4}, we have two pairs of adjoint functors 
$(K,R)$ and $(K',R')$. Explicitely
$$K:\ \Mm_B\to \Mm^{\Dd},~~K(N)=N\ot_B\Sigma;$$
$$R:\ \Mm^{\Dd}\to \Mm_B,~~R(M)=\Hom^{\Dd}(\Sigma, M)\cong M\ot^\Dd \Sigma^*.$$
The unit and counit will be called $\eta$ and $\varepsilon$, and are given
by the formulas
$$\eta_N:\ N\to \Hom^\Dd(\Sigma, N\ot_B\Sigma),~~
\eta_N(n)(u)=n\ot_B u;$$
$$\varepsilon_M:\ \Hom^\Dd(\Sigma,M)\ot_B\Sigma\to M,~~ 
\varepsilon_M(\varphi\ot_B u)=\varphi(u).$$
or
$$\eta_N:\ N\to (N\ot_B\Sigma)\ot^\Dd\Sigma^*,~~
\eta_N(n)=n\ot_Be;$$
$$\varepsilon_M:\ (M\ot^\Dd \Sigma^*)\ot_B\Sigma\to M,~~ 
\varepsilon_M(\sum_j m_j\ot_A g_j\ot_B u_j)=\sum_j m_jg_j(u_j).$$
The functors $K'$ and $R'$ are defined in a similar way.

For every $N\in \Mm_B$, we will consider the map 
$l_N=N\ot_B l:\ N\to N\ot_B\Sigma\ot_A
\Sigma^*$, $l_N(n)=n\ot_B e$.

\begin{definition}\delabel{2.1}
Let $B$ be a ring. We call $P\in {}_B\Mm$ totally faithful if for
all $N\in \Mm_B$ and $n\in N$, we have
\begin{equation}\eqlabel{2.1.1}
n\ot_B p=0~{\rm in}~N\ot_BP,~{\rm for~all}~p\in P~~\Longrightarrow~~
n=0.
\end{equation}
\end{definition}

Observe that $P$ is a faithful module if \equref{2.1.1} holds for
$N=B$; in fact total faithfulness is a purity condition.

\begin{lemma}\lelabel{2.2}
Let $\Sigma\in {}_B\Mm_A$ finitely generated projective as a right $A$-module. 
Then $\Sigma$ is totally faithful as a left
$B$-module if and only if $l:\ B\to \End_A(\Sigma)\cong
\Sigma\ot_A \Sigma^*$ is pure as a morphism of left $B$-modules.
\end{lemma}

\begin{proof}
Assume first that $\Sigma$ is totally faithful. Observe that
$l(b)=b\ot_Be=\sum_i be_i\ot_A f_i$. Take $N\in \Mm_B$, and $n\in N$. If
$$(I_N\ot_B l)(n\ot_B 1_B)=\sum_i n\ot_B e_i\ot_A f_i=0,$$
then for all $u\in \Sigma$, $0=\sum_i n\ot_B e_if_i(u)=n\ot_B u$, hence
$n=0$, and it follows that $I_N\ot_B l$ is injective, hence $l$ is pure.\\
Conversely, assume that $I_N\ot_B l$ is injective, for all $N\in \Mm_B$.
If $n\ot_B u=0$, for all $u\in \Sigma$, then $\sum_i n\ot_B e_i\ot_A f_i =0$,
hence $n=0$.
\end{proof}

\begin{proposition}\prlabel{2.3}
The functor $K$ is fully faithful if and only if $\Sigma$ is totally faithful
as a left $B$-module if and only if $l:\ B\to \Sigma\ot_A\Sigma^*$ is pure
in ${}_B\Mm$.
\end{proposition}

\begin{proof}
Take $N\in \Mm_B$.
The map $i_N=I_N\ot_B l:\ N\to N\ot_B\Sigma\ot_A\Sigma^*$ factorizes through
$\eta_N:\ N\to (N\ot_B\Sigma)\ot^\Dd\Sigma^*$. Hence $i_N$ is injective if
and only if $\eta_N$ is injective.\\
If $K$ is fully faithful, then every $\eta_N$ is bijective, hence injective,
hence every $i_N$ is injective, and $\Sigma$ is totally faithful.\\
Conversely, let $\Sigma\in {}_B\Mm$ be totally faithful, and take
$N\in \Mm_B$. We already know that $\eta_N$ is injective, and we are done
if we can show that it is also surjective. Consider
$\tilde{N}=(N\ot_B\Sigma)\ot^\Dd\Sigma^*/\eta_N(N)$,
and the canonical projection
$$\pi:\ (N\ot_B\Sigma)\ot^\Dd\Sigma^*\to \tilde{N}.$$
Let $x=\sum_j n_j\ot_B u_j\ot_A g_j\in( N\ot_B\Sigma)\ot^\Dd\Sigma^*$. Then
$$\sum_j \eta(n_j)\ot_B u_j\ot_A g_j=
\sum_{i,j} n_j\ot_B e_i\ot_A f_i\ot_B u_j\ot_A g_j=
\sum_{i} x\ot_B e_i\ot_A f_i.$$
Applying $\pi$ to the first three tensor factors, we find
$$0=\sum_j \pi(\eta(n_j))\ot_B u_j\ot_A g_j=
\sum_{i}\pi( x)\ot_B e_i\ot_A f_i,$$
hence for all $u\in \Sigma$,
$$0=\sum_{i}\pi( x)\ot_B e_if_i(u)=\sum_{i}\pi( x)\ot_Bu,$$
so $\pi(x)=0$, and $x\in \im(\eta_N)$, as needed.
\end{proof}

We now want to investigate when $R$ is fully faithful, or, equivalently,
when is $\varepsilon$ a natural isomorphism. For $M\in \Mm^\Dd$, we have
inclusions
$$(M\ot^\Dd\Sigma^*)\ot_B\Sigma\rTo^jM\ot^{\Dd}(\Sigma^*\ot_B\Sigma)
\subset M\ot_A\Sigma^*\ot_B\Sigma,$$
and an isomorphism
$$I_M\ot_A\varepsilon_\Dd:\ M\ot^\Dd(\Sigma^*\ot_B\Sigma)\to M.$$
It is obvious that $\varepsilon_M=(I_M\ot_A\varepsilon_\Dd)\circ j$,
hence $\varepsilon_M$ is an isomorphism if and only if $j$ is an 
isomorphism. Since $M\ot^\Dd \Sigma^*$ is the equalizer of
$\rho^r_M\ot_A \Sigma^*$ and $M\ot_A\rho^l_{\Sigma^*}= l_{M\ot_A\Sigma^*}$,
we have the following result.

\begin{proposition}\prlabel{2.4}
For $M\in \Mm^\Dd$, the following assertions are equivalent,
\begin{enumerate}
\item $j:\ (M\ot^\Dd\Sigma^*)\ot_B\Sigma\to M\ot^{\Dd}(\Sigma^*\ot_B\Sigma)$
is an isomorphism;
\item $\varepsilon_M$ is an isomorphism;
\item $\bullet\ot_B\Sigma$ preserves the equalizer of
$\rho^r_M\ot_A \Sigma^*$ and $l_{M\ot_A\Sigma^*}$.
\end{enumerate}
$R$ is fully faithful if and only if these three conditions are satisfied
for every $M\in \Mm^\Dd$. In particular, $R$ is fully faithful if
$\Sigma\in \Mm_B$ is flat.
\end{proposition}

We will now give some sufficient conditions for the fully
faithfulness of $R$. First we make the following general observation. Let $B$ be
a ring, $P$ and $Q$ $B$-bimodules, and $M$ a right $B$-module.
Then $\Hom_B(M,P)\in {}_B\Mm$, with left $B$-action $(b\cdot\alpha)(m)=
b\alpha(m)$, for all $m\in M$, $b\in B$ and right $B$-linear
$\alpha:\ M\to P$.
If $f:\ P\to Q$ is left and right $B$-linear, then
$$f\circ\bullet:\ \Hom_B(M,P)\to \Hom_B(M,Q)$$
is left $B$-linear. If $f$ is a split epimorphism in ${}_B\Mm_B$, split
by $g:\ Q\to P$, then $f\circ\bullet$ is a split epimorphism in ${}_B\Mm$,
split by $g\circ\bullet$.\\

We consider the contravariant functor
$C=\Hom_\ZZ(\bullet,\QQ/\ZZ):\ \Ab\to \Ab.$
$\QQ/\ZZ$ is an injective cogenerator of $\Ab$, and therefore
$C$ is exact and reflects isomorphisms. If $B$ is a ring, then
$C$ induces functors
$$C:\ \Mm_B\to {}_B\Mm~~{\rm and}~~{}_B\Mm\to \Mm_B.$$
For example, if $M\in \Mm_B$, then $C(M)$ is a left $B$-module,
by putting $(b\cdot f)(m)=f(mb)$.
For $M\in \Mm_B$ and $P\in{}_BM$, we have the following isomorphisms,
natural in $M$ and $P$:
\begin{equation}\eqlabel{2.5.1a}
\Hom_B(M,C(P))\cong {}_B\Hom(P,C(M))\cong C(M\ot_B P)
\end{equation}
If $P\in {}_B\Mm_B$, then $C(P)\in {}_B\Mm_B$, and the above isomorphisms
are isomorphisms of left $B$-modules.\\
The map $l:\ B\to \Sigma\ot_A \Sigma^*$ is a $B$-bimodule map. Hence the
map $C(l):\ C(\Sigma\ot_A \Sigma^*)\to C(B)$ is also a $B$-bimodule map.

\begin{proposition}\prlabel{2.5a}
Let $A$ and $B$ be rings, and $\Sigma\in {}_B\Mm_{A,{\rm fgp}}$.
If $C(l):\ C(\Sigma\ot_A \Sigma^*)\to C(B)$ is a split epimorphism in
${}_B\Mm_B$, then the functors $R$ and $R'$ are both fully faithful.
\end{proposition}

\begin{proof}
We will show that $R$ is fully faithful; the proof of the fact that $R'$ is
fully faithful is similar.
From \prref{2.4}, it follows that it suffices to show that the sequence
\begin{equation}\eqlabel{2.5.1}
0\to (M\ot^\Dd \Sigma^*)\ot_B\Sigma\to M\ot_A\Dd
\pile{\rTo^{\rho\ot_A \Sigma^*\ot_B \Sigma}\\ \rTo_{l_{M\ot_A\Sigma^*}\ot_B\Sigma}}
M\ot_A\Dd\ot_A\Dd
\end{equation}
is exact, for every $(M,\rho)\in \Mm^\Dd$.\\
If $C(l):\ C(\Sigma\ot_A \Sigma^*)\to C(B)$ is a split $B$-bimodule epimorphism,
then it follows from the remarks preceding the Proposition that
$$C(l)\circ\bullet:\ \Hom_B(M,C(\Sigma\ot_A\Sigma^*))\to \Hom_B(M,C(B))$$
is a split epimorphism in ${}_B\Mm$, for every $M\in \Mm_B$.
Applying \equref{2.5.1a}, we find that
$$C(l_M):\ C(M\ot_B\Sigma\ot_A\Sigma^*)\to C(M)$$
is a split epimorphism in ${}_B\Mm$. Now consider the  following diagram
in $\Mm_B$.
$$\begin{diagram}
0&&0\\
\dTo&&\dTo\\
M\ot^\Dd\Sigma^*&\rTo^{j}&M\ot_A\Sigma^*\\
\dTo^{j}&&\dTo^{\rho\ot_A\Sigma^*}\\
M\ot_A\Sigma^*&\rTo^{l_{M\ot_A\Sigma^*}}&M\ot_A\Dd\ot_A\Sigma^*\\
\dTo^{\rho\ot_A\Sigma^* }\dTo_{l_{M\ot_A\Sigma^*}}&&
\dTo^{\rho \ot_A\Dd\ot_A\Sigma^*}\dTo_{l_{M\ot_A\Sigma^*}\ot_B\Sigma\ot_A\Sigma^*}\\
M\ot_A\Dd\ot_A\Sigma^*&\rTo^ {l_{M\ot_A\Dd\ot_A\Sigma^*}}&
M\ot_A\Dd\ot_A\Dd\ot_A\Sigma^*
\end{diagram}$$
A straightforward computation shows that the two squares in the diagram
commute. It is also easy to see that the right column is exact: take
$$x=\sum_j m_j\ot_Ag_j\ot_B u_j\ot_A h_j\in M\ot_A\Dd\ot_A\Sigma^*,$$
and assume that $x$ lies in the equalizer of 
$\rho \ot_A\Dd\ot_A\Sigma^*$ and $l_{M\ot_A\Sigma^*}\ot_B\Sigma\ot_A\Sigma^*$.
Then
$$\sum_j \rho(m_j)\ot_Ag_j\ot_B u_j\ot_A h_j=
\sum_{j,i} m_j\ot_Ag_j\ot_B e_i\ot_A f_i\ot_B u_j\ot_A h_j,$$
hence
\begin{eqnarray*}
&&\hspace*{-2cm}
x= \sum_j m_j\ot_Ag_j\ot_B u_j\ot_A h_j
= \sum_{j,i} m_j\ot_Ag_j\ot_B e_if_i(u_j)\ot_A h_j\\
&=& \sum_j \rho(m_j)g_j( u_j)\ot_A h_j
= \sum_j \rho(m_j)\ot_A g_j( u_j) h_j\\
&=& (\rho\ot_A\Sigma^*)(\sum_j m_j\ot_A g_j( u_j) h_j).
\end{eqnarray*}
Now we apply the functor $C$ to the above diagram. Then we obtain a commutative
diagram in ${}_B\Mm$, with exact columns.
$$\begin{diagram}
C(M\ot_A\Dd\ot_A\Dd\ot_A\Sigma^*)&
\pile{\lTo^{h'}\\ \rTo_{C(l_{M\ot_A\Dd\ot_A\Sigma^*})}}&
C(M\ot_A\Dd\ot_A\Sigma^*)\\
\dTo^{C(\rho\ot_A\Dd\ot_A\Sigma^*)}\dTo_{C(l_{M\ot_A\Sigma^*}\ot_B\Sigma\ot_A\Sigma^*)}
&&
\dTo^{C(\rho\ot_A\Sigma^* )}\dTo_{C(l_{M\ot_A\Sigma^*})}\\
C(M\ot_A\Dd\ot_A\Sigma^*)&
\pile{\lTo^{h}\\ \rTo_{C(l_{M\ot_A\Sigma^*})}}&
C(M\ot_A\Sigma^*)\\
\dTo^{C(\rho\ot_A\Sigma^* )}&&\dTo_{C(j)}\\
C(M\ot_A\Sigma^*)&
\pile{\lTo^{k}\\ \rTo_{C(j)}}&
C(M\ot^\Cc\Sigma^*)\\
\dTo&&\dTo\\
0&&0
\end{diagram}$$
We know from the above arguments that $C(l_{M\ot_A\Dd\ot_A\Sigma^*})$
and $C(l_{M\ot_A\Sigma^*})$ have right inverses in ${}_B\Mm$. Diagram chasing
arguments then show that $C(j)$ has a right inverse $k$ in ${}_B\Mm$
such that $k\circ C(j)=C(\rho\ot_A\Sigma^*)\circ h$. Thus the bottom
row in the above diagram is a split fork, split by
$$C(M\ot_A\Dd\ot_A\Sigma^*)\lTo^{h}C(M\ot_A\Sigma^*)\lTo^{k}
C(M\ot^\Dd\Sigma^*)$$
(see \cite[p.149]{McLane} for the definition of a split fork).
Split forks are preserved by arbitrary functors, so applying
${}_B\Hom(A,\bullet)$, we obtain a split fork in ${}_B\Mm$. Using \equref{2.5.1},
this split fork takes the form
$$C(M\ot_A\Dd\ot_A\Dd)
\pile{\rTo^{C(\rho \ot_A\Sigma^*\ot_B\Sigma)}\\ \rTo_{C(l_{M\ot_A\Sigma^*}\ot_B\Sigma)}}
C(M\ot_A\Dd)\rTo^{C(j\ot_B\Sigma)}
C((M\ot^\Dd\Sigma^*)\ot_B\Sigma).$$
$C$ is exact and reflects isomorphisms, hence it also reflects coequalizers.
It then follows that \equref{2.5.1} is exact, and we are done.
\end{proof}

\begin{proposition}\prlabel{2.5}
Let $A$ and $B$ be rings, and $\Sigma\in {}_B\Mm_{A,{\rm fgp}}$. 
If $\Sigma\in {}_B\Mm$ is totally faithful, and $l$ maps $B$ into the center of
$\Sigma\ot_A\Sigma^*\cong \End_A(\Sigma)$, or, equivalently, every $\varphi
\in \End_A(\Sigma)$ is left $B$-linear,
then $C(l):\ C(\Sigma\ot_A \Sigma^*)\to C(B)$ is a split epimorphism in
${}_B\Mm_B$, and the functors $R$ and $R'$ are fully faithful.
\end{proposition}

\begin{proof}
From \leref{2.2}, we know that
$l:\ B\to \Sigma\ot_A\Sigma^*$, $l(b)=\sum_i e_i\ot_A f_ib=
\sum_i be_i\ot_A f_i$ is pure in ${}_B\Mm$. This means that, for every
$N\in\Mm_B$, the map 
$$l_N:\ N\to N\ot_B\Sigma\ot_A\Sigma^* ,~~l_N(n)=\sum_i n\ot_Be_i\ot_A f_i,$$
is injective. In particular, $l_{C(B)}$ is an injective left $B$-linear map.
Applying the contravariant functor $C$, we find that
$$C(l_{C(B)}):\ C(C(B)\ot_B\Sigma\ot_A\Sigma^*)\to C(C(B))$$
is an epimorphism in $\Mm_B$. From \equref{2.5.1a}, it then follows that
$$C(l)\circ\bullet:\
\Hom_B(C(B), C(\Sigma\ot_A\Sigma^*))\to \Hom_B(C(B), C(B))$$
is an epimorphism, which implies that $C(l)$ is a split epimorphism in
$\Mm_B$. The condition $l(B)\subset Z(\End_A(\Sigma))$ is equivalent to
$(\Sigma\ot_A \Sigma^*)^B=\Sigma\ot_A \Sigma^*$. Since $l$ is injective,
it also implies that $B$ is commutative, hence $B^B=B$. It follows that
$$C(B)^B=C(B)~~{\rm and}~~(C(\Sigma\ot_A \Sigma^*))^B=C(\Sigma\ot_A \Sigma^*),$$
so $C(l)$ is also a split epimorphism in ${}_B\Mm_B$.
\end{proof}

Combining Propositions \ref{pr:2.3} and \ref{pr:2.5}, we obtain the following
result.

\begin{theorem}\thlabel{2.6} Let $A$ and $B$ be rings,  $\Sigma\in {}_B\Mm_A$
finitely generated and projective as a right $A$-module, with finite dual
basis $e$, and  $\Dd=\Sigma^*\ot_B
\Sigma$. Consider the adjoint pairs $(K,R)$ and $(K',R')$ introduced above.
If every $\varphi\in \End_A(\Sigma)$ is left $B$-linear,
then the following assertions are equivalent:
\begin{enumerate}
\item $(K,R)$ is a pair of inverse equivalences;
\item $(K',R')$ is a pair of inverse equivalences;
\item $K$ is fully faithful;
\item $K'$ is fully faithful;
\item $l:\ B\to \Sigma\ot_A\Sigma^*$, $l(b)=be=eb$, is pure in ${}_B\Mm$;
\item $\Sigma\in{}_B\Mm$ is totally faithful;
\item $l:\ B\to \Sigma\ot_A\Sigma^*$, $l(b)=be=eb$, is pure in $\Mm_B$;
\item $\Sigma^*\in\Mm_B$ is totally faithful
\end{enumerate}
\end{theorem}

In the situation where $\Sigma=A$, and $A$ and $B$ commutative, we recover
the Joyal-Tierney Theorem. Our proof is an adaption of Mesablishvili's proof,
see \cite{Mesablishvili}. In \cite[Prop. 2.3]{Caenepeel03}, the case where
$\Sigma=A$, and $A$ not necessarily commutative is discussed. We remark that
the formulation of \cite[Prop. 2.3]{Caenepeel03} is incorrect:
we need the additional assumption that $A^B=B$.
This was pointed out to us by B. Mesablishvili.\\

The next result is due to El Kaoutit and G\'omez Torrecillas (see
\cite[Theorem 3.10]{Kaoutit}). Mesablishvili informed us that it is a
special case of Beck's Theorem. We give a short proof, for completeness sake.

\begin{theorem} {\bf (Faithfully flat descent)}\thlabel{2.7}
Let $A$ and $B$ be rings,  $\Sigma\in {}_B\Mm_A$
finitely generated and projective as a right $A$-module and flat as a left $B$-module.
Then $(K,R)$ is a pair of inverse equivalences if and only if
$\Sigma\in {}_B\Mm$ is faithfully flat.
\end{theorem}

\begin{proof}
First assume that $\Sigma\in {}_B\Mm$ is faithfully flat. For any $N\in \Mm_B$,
the map
$$f:\ N\ot_B\Sigma\to N\ot_B\Sigma\ot_A\Sigma^*\ot_B\Sigma,~~
f(n\ot_B u)=\sum_i n\ot_B e_i\ot_Af_i\ot_B u$$
is injective: if
$f(\sum_j n_j\ot_B u_j)= \sum_{i,j} n_j\ot_B e_i\ot_Af_i\ot_B u_j=0$,
then
$0=\sum_{i,j} n_j\ot_B e_if_i(u_j)=\sum_j n_j\ot_b u_j$.
Since $\Sigma$ is faithfully flat, it follows that
$$l_N:\ N\to N\ot_B\Sigma\ot_A\Sigma^*,~~l(n)=\sum_i l\ot_Be_i\ot_A f_i$$
is injective, and this means that $l$ is pure. It then follows from \prref{2.3}
that $K$ is fully faithful.\\
Conversely, let $E$ be a short sequence
in $\Mm_B$ such that $E\ot_B\Sigma$
is exact in $\Mm_A$. Applying the exact functor $R$ to $E\ot_B\Sigma$, and using the
fact that $\eta$ is an isomorphism, we find that $E$ is exact, and it
follows that $\Sigma\in {}_B\Mm$ is faithfully flat.
\end{proof}

\section{Galois corings}\selabel{3}
Let $A$ and $B$ be rings, $\Cc$ an $A$-coring, and $\Sigma\in {}_B\Mm^\Cc_{\rm fgp}$,
and consider the adjoint pair of functors $(F,G)$ introduced in \seref{1}. We can
then also consider the comatrix coring $\Dd=\Sigma^*\ot_B\Sigma$. We will now discuss
when $(F,G)$ is a pair of inverse equivalences.

\begin{lemma}\lelabel{3.1}
The map
$${\rm can}:\ \Dd\to\Cc,~~{\rm can}(g\ot_B u)=g(u_{[0]})u_{[1]}$$
is a morphism of corings.
\end{lemma}

\begin{proof}
It is obvious that ${\rm can}$ is an $A$-bimodule map. We also compute that
\begin{eqnarray*}
&&\hspace*{-10mm}
({\rm can}\ot_A{\rm can})(\Delta_\Dd(g\ot_B u))=
\sum_i {\rm can}(g\ot_B e_i)\ot_A {\rm can}(f_i\ot_B u)\\
&=& \sum_i g(e_{i[0]})e_{i[1]}\ot_A f_i(u_{[0]})u_{[1]}
= g(u_{[0]})u_{[1]}\ot_A u_{[2]}=\Delta_\Cc({\rm can}(g\ot_B u))
\end{eqnarray*}
and
$$\varepsilon_\Cc({\rm can}(g\ot_B u))=g(u)=\varepsilon_\Dd(g\ot_B u).$$
\end{proof}

\begin{lemma}\lelabel{3.1b}
We have a functor
$$\Gamma:\ \Mm^\Dd\to \Mm^\Cc,~~\Gamma(M,\tilde{\rho})=
(M,\rho= (M\ot_A{\rm can})\circ \tilde{\rho}).$$
$\Gamma\circ K= F$,
and we have a natural inclusion
$\alpha:\ R\to  G\circ\Gamma$.
If ${\rm can}$ is bijective, then $\Gamma$ is an isomorphism of categories,
and $\alpha$ is a natural isomorphism.
\end{lemma}

\begin{proof}
We know that $(\Sigma,\rho)\in \Mm^\Cc$, and $(\Sigma,\tilde{\rho})\in \Mm^\Dd$,
with $\tilde{\rho}(u)=\sum_i u\ot_A e_i\ot_b f_i$. We write $\rho(u)=
u_{[0]}\ot_A u_{[1]}$. Then $\Gamma(\Sigma,\tilde{\rho})= (\Sigma,\rho)$,
since
$$(M\ot_A{\rm can})(\tilde{\rho}(u))=
\sum_i e_i\ot_A f_i(u_{[0]})u_{[1]}=u_{[0]}\ot_A u_{[1]}.$$
Consequently $\Gamma(K(N))=\Gamma(N\ot_B \Sigma)=F(N)$,
for all $N\in \Mm_B$.
Now take $M\in \Mm^\Dd$ and
$f\in R(M)=\Hom^{\Dd}(\Sigma, M)$. Then $\Gamma(f)=f:\ \Gamma(\Sigma)\to \Gamma(M)$
is right $\Cc$-colinear, since for all $u\in \sigma$:
\begin{eqnarray*}
&&\hspace*{-10mm}
(f\ot_A\Cc)(\rho(u))=
\Bigl((f\ot_A\Cc)\circ (\Sigma\ot_A{\rm can})\circ \tilde{\rho}\Bigr)(u)\\
&=& \bigl((M\ot_A {\rm can})\circ (f\ot_A\Dd)\circ \tilde{\rho}\Bigr)(u)
= (M\ot_A {\rm can})(\tilde{\rho}(f(u))=\rho(f(u)),
\end{eqnarray*}
and we find that 
$$R(M)=\Hom^{\Dd}(\Sigma, M)\subset G(\Gamma(M))=\Hom^\Cc(\Gamma(M),\Gamma(\Sigma)).$$
The rest of the proof is obvious.
\end{proof}

As an immediate consequence, we have:

\begin{proposition}\prlabel{3.2}
With notation as above, if ${\rm can}$ is an isomorphism, then $F$ is fully
faithful if and only if $K$ is fully faithful, and $G$ is fully faithful
if and only if $R$ is fully faithful.
\end{proposition}

We now give some necessary conditions for $(F,G)$ to be a pair of
inverse equivalences.

\begin{proposition}\prlabel{3.3}
With notation as above, we have the following results.
\begin{enumerate}
\item If the functor $F$ is fully faithful, then the map $l:\ B\to T=\Sigma\ot^\Cc\Sigma^*$,
$l(b)=eb=be$ is an isomorphism;
\item if the functor $G$ is fully faithful, then the map ${\rm can}:\ \Dd\to\Cc$
is an isomorphism.
\end{enumerate}
\end{proposition}

\begin{proof}
1) This follows from the observation that $l=\nu_B$.\\
2) From \leref{1.2b}, we have an isomorphism $\alpha:\ \Sigma^*\to \Hom^\Cc(\Sigma,\Cc)$.
We easily check that
${\rm can}=\zeta_\Cc\circ (\alpha\ot_B \Sigma)$,
hence ${\rm can}$ is an isomorphism if and only if $\zeta_\Cc$ is an isomorphism.
\end{proof}

\begin{definition}\delabel{3.4} (\cite[3.4]{Kaoutit})
Let $\Cc$ be an $A$-coring, $\Sigma\in \Mm^\Cc_{\rm fgp}$, and let
$T=\Sigma\ot^\Cc \Sigma^*\cong \End^\Cc(\Sigma)$. Then we call $(\Cc,\Sigma)$
a Galois coring if ${\rm can}:\ \Dd=\Sigma^*\ot_T \Sigma\to \Cc$ is an isomorphism.
\end{definition}

Let us remark that a diffferent terminology is used in \cite{Br3}. If $(\Cc,\Sigma)$ is
a Galois coring in the sense of \deref{3.4}, then $\Sigma$ is called a Galois $\Cc$-comodule.

We will now give some equivalent definitions. Recall first that $(M,\rho)\in \Mm^\Cc$
is termed $(\Cc,A)$-injective if the following holds: for every
right $\Cc$-colinear map $i:\ N\to L$ having a left inverse in $\Mm_A$,
and for every $f:\ N\to M$ in $\Mm^\Cc$, there exists a $g:\ L\to M$ in
$\Mm^\Cc$ such that $g\circ i= f$. An easy computation shows that $(M,\rho)$ is
$(\Cc,A)$-injective if and only if $\rho$ has a left inverse in $\Mm^\Cc$.

\begin{proposition}\prlabel{3.5}
Let $\Cc$ be an $A$-coring and $\Sigma\in \Mm^\Cc_{\rm fgp}$. Then the
following assertions are equivalent.
\begin{enumerate}
\item $(\Cc,\Sigma)$ is Galois;
\item the evaluation map
${\rm ev}_\Cc:\ \Hom^\Cc(\Sigma,\Cc)\ot_T \Sigma\to \Cc$
is an isomorphism;
\item if $M\in \Mm^\Cc$ is $(\Cc,A)$-injective, then the evaluation map
$${\rm ev}_M:\ \Hom^\Cc(\Sigma,M)\ot_T \Sigma\to M,~~
{\rm ev}_M(f\ot_T u)=f(u)$$
is an isomorphism.
\end{enumerate}
\end{proposition}

\begin{proof}
$1)\Longleftrightarrow 2)$ follows from the fact that
$\Hom^\Cc(\Sigma,\Cc)\cong \Hom_A(\Sigma, A)=\Sigma^*$, see \leref{1.2b}.
$3)\Longrightarrow 2)$ is obvious.\\
$1)\Longrightarrow 3)$. For all $L\in \Mm^\Cc$, we have a split exact
sequence (see \cite[3.7]{W3}):
\begin{equation}\eqlabel{3.5.1}
0\to \Hom^\Cc(L,M)\rTo^{i}\Hom_A(L,M)\rTo^{j}\Hom_A(L,M\ot_A\Cc).
\end{equation}
The map $j$ is given by
$$j(f)(l)=f(l)_{[0]}\ot_A f(l)_{[1]}-f(l_{[0]})\ot_A l_{[1]},$$
and the splitting maps $\alpha:\ \Hom_A(L,M)\to \Hom^\Cc(L,M)$
and $\beta:\ \Hom_A(L,M\ot_A\Cc)\to \Hom_A(L,M)$ are given by the formulas
$$\alpha(f)(l)=\gamma(f(l_{[0]})\ot_A l_{[1]}),~~{\rm and}~~
\beta(g)=\gamma\circ g,$$
where $\gamma$ is a left inverse of $\rho$ in $\Mm^\Cc$. Now take $L=\Sigma$,
and apply $\bullet_A \Sigma$ to \equref{3.5.1}. Using the fact that
$\Hom_A(\Sigma,M)\cong M\ot_A\Sigma^*$ and $\Hom_A(\Sigma,M\ot_A\Cc)\cong
M\ot_A\Cc\ot_A\Sigma^*$, we obtain a diagram
$$\begin{diagram}
0&\to&\Hom^{\Cc}(\Sigma,M)\ot_T\Sigma&\to&M\ot_A\Dd&\to& M\ot_A\Cc\ot_A \Dd\\
&&\dTo_{{\rm ev}_M}&&\dTo_{M\ot_A{\rm can}}&&\dTo_{M\ot_A\Cc\ot_A{\rm can}}\\
0&\to& M&\rTo^{\rho}&M\ot_A\Cc&\rTo^{\psi}&M\ot_A\Cc\ot_A \Cc
\end{diagram}$$
where $\psi=\rho\ot_A \Cc-M\ot_A\Delta_\Cc$. The toprow is split exact,
and the bottomrow is exact. a straightforward computation shows that
the diagram commutes. From the fact that ${\rm can}$ is bijective, it then
follows that ${\rm ev}_M$ is bijective.
\end{proof}

We now look at corings with a fixed flat comodule. First we have
to recall basic facts about generators. We include the proof of our
next Lemma for completeness sake.

\begin{lemma}\lelabel{3.7}
Let $\Cc$ be an $A$-coring, and $\Sigma\in \Mm^\Cc$, and consider the following
statements.
\begin{enumerate}
\item $\Sigma$ generates $\Mm^\Cc$: if $0\neq g:\ M\to N$ in $\Mm^\Cc$,
then there exists $f\in \Hom^\Cc(\Sigma,M)$ such that $g\circ f\neq 0$;
\item for all $M\in \Mm^\Cc$, ${\rm ev}_M:\ \Hom^\Cc(\Sigma,M)\ot_B \Sigma
\to M$ is surjective;
\item for all $M\in \Mm^\Cc$, ${\rm ev}_M:\ \Hom^\Cc(\Sigma,M)\ot_B \Sigma
\to M$ is bijective.
\end{enumerate}
The first two statements are equivalent. If $\Cc$ is flat as a left $A$-module,
then all three statements are equivalent.
\end{lemma}

\begin{proof}
$\ul{1)\Rightarrow 2)}$. The image of ${\rm ev}_M$ is a right $\Cc$-comodule,
and we can consider the canonical projection $g:\ M\to M/\im({\rm ev}_M)$ in $\Mm^\Cc$.
For all $f\in \Hom^\Cc(\Sigma,M)$ and $u\in \Sigma$,
$(g\circ f)(u)=g({\rm ev}_M(f\ot u))=0$, hence $g=0$, and ${\rm ev}_M$
is surjective.\\
$\ul{2)\Rightarrow 1)}$.
Take $m\in M$ such that $g(m)\neq 0$. We have $f_i\in \Hom^\Cc(\Sigma,M)$
and $u_i\in \Sigma$ such that $m=\sum_i f_i(u_i)$. If $g(\varphi_i(u_i))=0$
for all $i$, then $g(m)=g(\sum_i f_i(u_i))=0$, which is impossible.
Hence there exists $i$ such that $g\circ f_i\neq 0$.\\
$\ul{2)\Rightarrow 3)}$. (along the lines of \cite[43.12]{BrzezinskiWisbauer}).
Assume that $\Cc$ is flat as left $A$-module. We have to show that
every ${\rm ev}_M$ is injective. Take $\sum_{i=1}^k f_i\ot m_i\in
\Ker {\rm ev}_M$, i.e. $\sum_{i=1}^k f_i(m_i)=0$.
Consider the projection $\pi_i:\ \Sigma^k\to \Sigma$ onto the $i$-th component,
and
$f=\sum_{i=1}^k f_i\circ \pi_i\in \Hom^\Cc(\Sigma^k,M)$.
$\Ker f\in \Mm^\Cc$, since $\Cc$ is flat (see \cite{BrzezinskiWisbauer}).
Also $(m_1,\cdots, m_k)\in \Ker f$, since $f(x_1,\cdots,x_k)=
\sum_{i=1}^k f_i(x_i)$. By assumption, the map
$${\rm ev}_{\Ker f}:\ \Hom^\Cc(\Sigma,\Ker f)\ot_B\Sigma\to \Ker f$$
is surjective, hence we can find $a_j\in \Sigma$ and
$g_j\in \Hom^\Cc(\Sigma,\Ker f)$ such that
$\sum_{j=1}^l g_j(a_j)=(m_1,\cdots,m_k)$. Using the fact that
$\im g_j\subset \ker f$, we obtain
\begin{eqnarray*}
&&\hspace*{-2cm}
\sum_{i=1}^k f_i\ot_B m_i=\sum_{i=1}^k f_i\ot_B \sum_{j=1}^l(\pi\circ g_j)(a_j)\\
&=& \sum_{j=1}^l\left(\sum_{i=1}^k f_i\circ \pi_i\right)\circ g_j\ot_B a_j
= \sum_{j=1}^l f\circ g_j\ot_B a_j=0.
\end{eqnarray*}
\end{proof}

\begin{theorem}\thlabel{3.8}
Let $\Cc$ be an $A$-coring, $\Sigma\in \Mm^\Cc_{\rm fgp}$, and
$B=T=\End^\Cc(\Sigma)$. The following assertions are equivalent.
\begin{enumerate}
\item $(\Cc,\Sigma)$ is Galois and $\Sigma\in {}_B\Mm$ is flat;
\item $G$ is fully faithful and $\Sigma\in {}_B\Mm$ is flat;
\item $\Sigma\in \Mm^\Cc$ is a generator and $\Cc\in {}_A\Mm$ is flat.
\item ${\rm ev}_M$ is bijective for every $M\in \Mm^\Cc$ and
$\Sigma\in {}_B\Mm$ is flat.
\end{enumerate}
\end{theorem}

\begin{proof}
$\ul{1) \Rightarrow 2)}$ follows from Propositions \ref{pr:2.4} and
\ref{pr:3.2}.\\
$\ul{2) \Rightarrow 1)}$ follows from \prref{3.3}.\\
$\ul{2) \Rightarrow 3)}$. $\Sigma\in {}_B\Mm$ is flat, and
$\Sigma^*\in {}_A\Mm$ is finitely generated projective, hence flat,
so $\Sigma^*\ot_B\Sigma=\Dd\cong \Cc$ is flat in ${}_A\Mm$.\\
Take $0\neq g:\ M\to N$ in $\Mm^\Cc$. Then 
$$G(g):\ \Hom^\Cc(\Sigma,M)\to \Hom^\Cc(\Sigma,N),~~
G(g)(f)=g\circ f.$$
$G(g)\neq 0$ since $G$ is fully faithful. Hence there exists $f\in
\Hom^\Cc(\Sigma,M)$ such that $G(g)(f)=g\circ f\neq 0$, and this is
exactly what we need.\\
$\ul{3) \Rightarrow 4)}$ (along the lines of \cite[15.9]{W1}).
We first show that $\Sigma$ is flat as a left
$B$-module. It suffices to show (cf. e.g. \cite[12.16]{W1}) that, for any finitely generated
right ideal $J=f_1B+\cdots f_kB$ of $B$, the map
$\mu_J:\ J\ot_B\Sigma\to J\Sigma,~\mu_J(g\ot u)=g(u)$
is injective. We consider the surjection
$$\phi:\ \Sigma^n\to J\Sigma,~~\phi(u_1,\cdots,u_n)=\sum_{i=1}^n f_i(u_i)$$
$K=\Ker\phi\in \Mm^\Cc$, because $\Cc\in {}_A\Mm$ is flat. We have an
exact sequence
$$0\to \Hom^\Cc(\Sigma,K)\rTo^{\alpha} \Hom^\Cc(\Sigma,\Sigma^n)
\rTo^\beta \Hom^\Cc(\Sigma,J\Sigma)\to 0$$
$\alpha$ is the natural embedding, and $\beta(f)=\phi\circ f$. 
Observe that $\Hom^\Cc(\Sigma,J\Sigma)\cong J$. Tensoring by $\Sigma$,
we obtain the
following commutative diagram with exact rows:
$$\begin{diagram}
&&\Hom^\Cc(\Sigma,K)\ot_B\Sigma&\rTo^{\alpha\ot \Sigma}&\Hom^\Cc(\Sigma,\Sigma^n)
\ot_B\Sigma&
\rTo^{\beta\ot\Sigma}& J\ot_B\Sigma&\rTo& 0\\
&&\dTo_{{\rm ev}_K}&&\dTo_{{\rm ev}_{\Sigma^n}}&&\dTo_{\mu_J}&&\\
0&\rTo&K&\rTo &\sigma^n&\rTo^{\phi}&J\Sigma&\rTo& 0
\end{diagram}$$
${{\rm ev}_K}$ is surjective, by assumption, and ${\rm ev}_{\Sigma^n}$
is the canonical isomorphism 
$$\Hom^\Cc(\Sigma,\Sigma^n)\ot_B \Sigma\cong \Hom^\Cc(\Sigma,\Sigma)^n\ot_B \Sigma\cong 
B^n\ot_B\Sigma\cong \Sigma^n$$
A diagram chasing argument then implies that $\mu_J$ is injective.
It then follows from \leref{3.7} that every ${\rm ev}_M$ is bijective.\\
$\ul{4) \Rightarrow 1)}$ follows from \prref{3.5}.
\end{proof}

\begin{theorem}\thlabel{3.9}
Let $\Cc$ be an $A$-coring, $\Sigma\in \Mm^\Cc_{\rm fgp}$, and
$B=T=\End^\Cc(\Sigma)$. The following assertions are equivalent.
\begin{enumerate}
\item $(\Cc,\Sigma)$ is Galois and $\Sigma\in {}_B\Mm$ is faithfully flat;
\item $(F,G)$ is a pair of inverse equivalences and $\Sigma\in {}_B\Mm$ is flat;
\item $\Sigma\in \Mm^\Cc$ is a progenerator and $\Cc\in {}_A\Mm$ is flat.
\end{enumerate}
\end{theorem}

\begin{proof}
$\ul{1) \Rightarrow 2)}$. $\Sigma\in {}_B\Mm$ is faithfully flat, hence
$(K,R)$ is a pair of inverse equivalences, by \thref{2.7}. 
It then follows from \prref{3.2} that $(F,G)$ is a pair of inverse
equivalences.\\
$\ul{2) \Rightarrow 1)}$. It follows from \prref{3.3} that $(\Cc,\Sigma)$
is Galois.\\
$\ul{1) \Rightarrow 3)}$. In view of \thref{3.8}, we only have to show
that $\Sigma\in \Mm^\Cc$ is projective. Take an epimorphism $f:\ M\to N$
in $\Mm^\Cc$. We know from \thref{3.8} that ${\rm ev}_M$ and ${\rm ev}_N$
are isomorphisms. We also have a commutative diagram
$$\begin{diagram}
\Hom^\Cc(\Sigma,M)\ot_B\Sigma&\rTo^{\Hom^\Cc(\Sigma,f)\ot_B\Sigma}&
\Hom^\Cc(\Sigma,N)\ot_B\Sigma\\
\dTo^{{\rm ev}_M}&&\dTo_{{\rm ev}_N}\\
M&\rTo^{f}&N
\end{diagram}$$
so it follows that $\Hom^\Cc(\Sigma,f)\ot_B\Sigma$ is surjective. From the
fact that $\Sigma$ is a faithfully flat left $B$-module, it then follows that
$\Hom^\Cc(\Sigma,f)$ is projective, hence $\Sigma$ is a projective object
in $\Mm^\Cc$.\\
$\ul{3) \Rightarrow 1)}$. It follows from \thref{3.8} that $(\Cc,\Sigma)$ is Galois
and that $\Sigma\in {}_B\Mm$ is faithfully flat. Arguments similar to the ones
in \cite[18.4 (3)]{W1} show that for any right ideal $J$ of $B$, the inclusion
$J\subset \Hom^\Cc(\Sigma,J\Sigma)$ is an equality. Details
are as follows. Take $g\in  \Hom^\Cc(\Sigma,J\Sigma)$.  Let
$\{u_1,\cdots,u_k\}$ be a set of generators of $\Sigma\in \Mm_A$, and
write $g(u_i)=f_i(u_i)$, with $f_i\in J$. Let $J'$ be the subideal of
$J$ generated by $\{f_1,\cdots,f_k\}$. Since $\{f_1(u_1),\cdots,f_k(u_k)\}$
generate $\im(g)$ as a right $A$-module, we have that $\im(g)\subset
J'\Sigma$. Let $\pi_i:\ M^k\to M$ and $e_i:\ M\to M^k$ be the natural
projection and inclusion. The map 
$f=\sum_{i=1}^k f_i\circ\pi_i:\ \Sigma^k\to J'\Sigma$
is surjective; since $\Sigma\in \Mm^\Cc$ is projective, there exists
$h=\sum_{j=1}^k e_j\circ h_j:\ \Sigma\to \Sigma^k$
such that 
$$g=f\circ h=\sum_{i,j=1}^k f_i\circ\pi_i\circ e_j\circ h_j=\sum_{i=1}^k f_i\circ h_i\in
J'\Sigma\subset J\Sigma$$
If $J\neq B$, then $\Hom^\Cc(\Sigma,J\Sigma)\neq \Hom^\Cc(\Sigma,\Sigma)$,
hence $J\Sigma\neq \Sigma$, and this proves that $\Sigma\in {}_B\Mm$
is faithfully flat using \cite[12.17]{W1}
\end{proof}

\begin{remark}
A more general version of \thref{3.9}, with a different proof, was given by
El Kaoutit and G\'omez Torrecillas in \cite[Theorem 3.2]{Kaoutit}. In particular,
Condition (3) of \thref{3.9} implies that $\Sigma\in \Mm_A$ is finitely generated
and projective.
\end{remark}

\section{Morita theory}\selabel{4}
\subsection*{The dual of the canonical map}
As before, let $A$ and $B$ be rings, $\Cc$ an $A$-coring,
and $\Sigma\in {}_B\Mm^\Cc_{\rm fgp}$. Let $T=\End^\Cc(\Sigma)$
and $\Dd=\Sigma^*\ot_B\Sigma$. Then ${}^*\Cc={}_A\Hom(\Cc,A)$ is
a ring, with multiplication defined by
$$(f\#g)(c)=g(c_{(1)}f(c_{(2)})).$$
In a similar way, $\Cc^*=\Hom_A(\Cc,A)$ is
a ring, with multiplication defined by
$$(f\#g)(c)=f(g(c_{(1)})c_{(2)}).$$
We have a ring isomorphism
$$\alpha:\ {}^*\Dd={}_A\Hom(\Sigma^*\ot_B\Sigma,A)\to {}_B\End(\Sigma)^{\rm op}$$
given by
$$\alpha(\varphi)(u)=\sum_i e_i\varphi(f_i\ot_B u)~~{\rm and}~~
\alpha^{-1}(\psi)(f\ot_B u)=f(\psi(u))$$
In a similar way, we have a ring isomorphism
$$\beta:\ \Dd^*=\Hom_A(\Sigma^*\ot_B\Sigma,A)\to \End_B(\Sigma^*)$$
given by
$$(\beta(\varphi)(f))(u)=\varphi(f\ot_B u)~~{\rm and}~~
\alpha^{-1}(\psi)(f\ot_B u)=\psi(f)(u)$$
Observe also that 
$${}_B\End(\Sigma)^{\rm op}\cong \End_B(\Sigma^*),$$
the isomorphism is given by sending $\psi$ to $\psi^*$.\\
We can also consider the maps dual to ${\rm can}:\ \Dd\to \Cc$:
$${}^*{\rm can}:\ {}^*\Cc\to {}_B\End(\Sigma)^{\rm op},~~
{}^*{\rm can}(\varphi(u)=u_{[0]}\varphi(u_{[1]})$$
$${\rm can}^*:\ \Cc^*\to \End_B(\Sigma^*),~~
{\rm can}^*(\varphi)(f)=\varphi\circ (f\ot \Cc)\circ \rho_{\Sigma}$$
We immediately have the following result:

\begin{proposition}\prlabel{4.1}
If $(\Cc,\Sigma)$ is Galois, then ${}^*{\rm can}$ and ${\rm can}^*$ are isomorphisms.\\
If ${}^*{\rm can}$ (resp. ${\rm can}^*$)
 is an isomorphism, and $\Cc\in {}_A\Mm$ (resp. $\Cc\in \Mm_A$)
 and $\Sigma\in {}_B\Mm$
are finitely generated projective, then $(\Cc,\Sigma)$ is Galois.
\end{proposition}

\subsection*{A Morita context associated to a comodule}
Let $\Cc$ be an $A$-coring, and $M\in {}^\Cc\Mm$. We can associate a Morita context
to $M$. If $\Cc=A$ is the trivial coring, then we recover the Morita context associated
to a module (see \cite{Ba}). The context will also generalize the Morita contexts
introduced in \cite{Abu} and \cite{CVW}.  The context will connect
$T={}^\Cc\End(M)^{\rm op}$ and ${}^*\Cc$.

\begin{lemma}\lelabel{4.2}
With notation as above, ${}^*M\in {}_T\Mm_{{}^*\Cc}$ and 
$Q={}^\Cc\Hom(\Cc,M)\in {}_{{}^*\Cc}\Mm_T$.
\end{lemma}

\begin{proof}
Let $\varphi\in {}^*M$, $f\in {}^*\Cc$, $t\in T$, $q\in Q$ and $m\in M$. The bimodule
structure on ${}^*M$ is defined by
\begin{equation}\eqlabel{4.2.1}
(\varphi\cdot f)(m)=f(m_{[-1]}\varphi(m_{[0]}))~~{\rm and}~~t\cdot\varphi=\varphi\circ t.
\end{equation}
Let us show that the two actions commute
\begin{eqnarray*}
&&\hspace*{-2cm}
(t\cdot (\varphi\cdot f))(m)=(\varphi\cdot f)(t(m))=
f\Bigl(t(m)_{[-1]}\varphi( t(m)_{[0]})\Bigr)\\
&=&f\Bigl(m_{[-1]}\varphi( t(m_{[0]}))\Bigr)=((t\cdot \varphi)\cdot f)(m).
\end{eqnarray*}
The bimodule
structure on $Q$ is defined by
\begin{equation}\eqlabel{4.2.2}
(f\cdot q)(c)=q(c_{(1)}f(c_{(2)}))~~{\rm and}~~q\cdot t=t\circ q.
\end{equation}
The two actions commute, since
$$((f\cdot q)\cdot t)(c)=t(q(c_{(1)}f(c_{(2)}))=(q\cdot t)(c_{(1)}f(c_{(2)}))=
(f\cdot (q\cdot t))(c).$$
\end{proof}

\begin{lemma}\lelabel{4.3}
With notation as in \leref{4.2}, we have well-defined bimodule maps
$$\mu:\ Q\ot_T {}^*M\to {}^*\Cc,~~\mu(q\ot\varphi)=\varphi\circ q;$$
$$\tau:\ {}^*M\ot_{{}^*\Cc}Q\to T,~~\tau(\varphi\ot q)(m)=q(m_{[-1]}\varphi(m_{[0]})).$$
\end{lemma}

\begin{proof}
These are straightforward verifications.
\end{proof}

\begin{theorem}\thlabel{4.4}
With notation as in Lemmas \ref{le:4.2} and \ref{le:4.3}, we have a Morita context
$\CC=(T, {}^*\Cc, {}^*M, Q,\tau,\mu)$.
\end{theorem}

\begin{proof}
We first show that $\mu\ot Q=Q\ot \tau$. For all $p, q\in Q$, $\varphi\in {}^*M$
and $c\in \Cc$,
we have
\begin{eqnarray*}
&&\hspace*{-2cm}
\Bigl((Q\ot\tau)(q\ot\varphi\ot p)\Bigr)(c)=(q\cdot \tau(\varphi\ot p))(c)=
\tau(\varphi\ot p)(q(c))\\
&=&p(q(c)_{[-1]}\varphi(q(c)_{[0]}))
= p(c_{(1)}\varphi(q(c_{(2)}))\\
&=&((\varphi\circ q)\cdot p)(c)=
\Bigl((\mu\ot q)(q\ot\varphi\ot p)\Bigr)(c)
\end{eqnarray*}
${}^*M\ot \mu=\tau\ot{}^*M$ since
\begin{eqnarray*}
&&\hspace*{-2cm}
\Bigl((\tau\ot{}^*M)(\varphi\ot q\ot\psi)\Bigr)(m)=(\tau(\varphi\ot q)\cdot\psi)(m)=
\psi(\tau(\varphi\ot q)(m))\\
&=&\psi(q(m_{[-1]}\varphi(m_{[0]}))=(\psi\circ q)(m_{[-1]}\varphi(m_{[0]}))\\
&=&(\varphi(\psi\circ q))(m)=
\Bigl(({}^*M\ot \mu)(\varphi\ot q\ot\psi)\Bigr)(m)
\end{eqnarray*}
for all $q\in Q$, $\varphi,\psi\in {}^*M$
and $c\in \Cc$.
\end{proof}

\begin{remarks}\relabel{4.5}
1) If $\Cc=A$ is the trivial coring, and $M\in {}_A\Mm$, then the Morita context
$\CC=({}_A\End(M)^{\rm op},A,{}^*M,M,\tau,\mu)$ is the Morita context associated
to the $A$-module $M$, as in \cite[II.4]{Ba}.\\
2) We can also construct a Morita context associated to $M\in \Mm^\Cc$:
$$\CC=(T=\End^\Cc(M),\Cc^*,Q=\Hom^\Cc(\Cc,M),M^*,\tau,\mu)$$
with $M^*\in {}_{\Cc^*}\Mm_T$ via
$$(f\cdot \varphi)(m)=f(\varphi(m_{[0]})m_{[1]})~~{\rm and}~~
\varphi\cdot t=\varphi\circ t,$$
$Q\in {}_T\Mm_{\Cc^*}$ via
$$(q\cdot f)(c)=q(f(c_{(1)})c_{(2)})~~{\rm and}~~t\cdot q=t\circ q.$$
The connecting maps are
$$\mu:\ M^*\ot_TQ\to \Cc^*,~~\mu(\varphi\ot q)=\varphi\circ q$$
$$\tau:\ Q\ot_{\Cc^*}M^*\to T,~~\tau(q\ot\varphi)(m)=q(\varphi(m_{[0]})m_{[1]})$$
3) Let $\Sigma\in \Mm^\Cc_{\rm fgp}$. Then $\Sigma^*\in {}^\Cc_{\rm fgp}\Mm$.
As ${}^\Cc\End(\Sigma^*)^{\rm op}\cong \End^\Cc(\Sigma)=T$, we obtain a
Morita context
\begin{equation}\eqlabel{4.5.1}
\CC=(T=\End^\Cc(\Sigma), {}^*\Cc,\Sigma,Q={}^\Cc\Hom(\Cc,\Sigma^*),
\tau,\mu)
\end{equation}
with $Q\in {}_{{}^*\Cc}\Mm_T$ by
$$(f\cdot q)(c)=q(c_{(1)}f(c_{(2)}))~~{\rm and}~~(q\cdot t)(c)=q(c)\circ t$$
and
$\Sigma\in {}_T\Mm_{{}^*\Cc}$ by
$$t\cdot u=t(u)~~{\rm and}~~u\cdot f=u_{[0]}f(u_{[1]})$$
and
\begin{equation}\eqlabel{4.5.2}
\mu:\ Q\ot_T\Sigma\to{}^*\Cc,~~\mu(q\ot u)(c)=q(c)(u)
\end{equation}
\begin{equation}\eqlabel{4.5.3}
\tau:\ \Sigma\ot_{{}^*\Cc} Q\to T,~~\tau(u\ot q)(v)=u_{[0]}(q(u_{[1]})(v))
\end{equation}
4) Take $x\in G(\Cc)$; then $A$ is a right $\Cc$-comodule: $\rho(a)=xa$.
The Morita context \equref{4.5.1} is then the Morita context studied in
\cite{Abu,CVW}.
\end{remarks}

If $\Sigma\in \Mm^\Cc$, then $\Sigma$ is also a right ${}^*\Cc$-module, and we
can associate to $\Sigma$ a Morita context as in \cite[II.4]{Ba}, namely
\begin{equation}\eqlabel{4.5.4}
\TT=(\widetilde{T}=\End_{{}^*\Cc}(\Sigma),{}^*\Cc,\Sigma,\widetilde{Q}=
\Hom_{{}^*\Cc}(\Sigma,{}^*\Cc),\widetilde{\tau},\widetilde{\varphi}),
\end{equation}
with
$$\widetilde{\mu}:\ \widetilde{Q}\ot_{\widetilde{T}}\Sigma\to {}^*\Cc,~~
\widetilde{\mu}(\lambda\ot u)=\lambda(u)$$
$$
\widetilde{\tau}:\ \Sigma\ot_{{}^*\Cc} \widetilde{Q}\to \widetilde{T},~~
\widetilde{\tau}(u\ot\lambda)(v)=u\cdot\lambda(v)=u_{[0]}(\lambda(v)(u_{[1]})).$$
We will now study the relationship between the Morita contexts \equref{4.5.1}
and \equref{4.5.4}. But first we need a Lemma.

\begin{lemma}\lelabel{4.6}
Consider $Q={}^\Cc\Hom(\Cc,\Sigma^*)$. A left $A$-linear map $q:\ \Cc\to \Sigma^*$
belongs to $Q$ if and only if
\begin{equation}\eqlabel{4.6.1}
c_{(1)}(q(c_{(2)})(u))=(q(c)(u_{[0]}))u_{[1]}
\end{equation}
for all $u\in \Sigma$.
\end{lemma}

\begin{proof} 
Recall that the left $\Cc$-coaction on $\Sigma^*$ is given by \equref{1.1.1}.
Hence $q\in Q$ if and only if
\begin{equation}\eqlabel{4.6.2}
c_{(1)}\ot_A q(c_{(2)})=\sum_i (q(c)(e_{i[0]}))e_{i[0]}\ot_A f_i.
\end{equation}
Applying the second tensor factor to $u\in \Sigma$, we obtain \equref{4.6.1}.
Conversely, if \equref{4.6.1} holds, then
\begin{eqnarray*}
&&\hspace*{-2cm} \sum_i (q(c)(e_{i[0]}))e_{i[1]}\ot_A f_i=
\sum_i c_{(1)}(q(c_{(2)})(e_i))\ot_A f_i\\
&=& \sum_i c_{(1)}\ot_A (q(c_{(2)})(e_i)) f_i=c_{(1)}\ot_A q(c_{(2)})
\end{eqnarray*}
proving \equref{4.6.2}.
\end{proof}

\begin{proposition}\prlabel{4.8}
We have a morphism of Morita contexts
$$\CC=(T,{}^*\Cc,\Sigma,Q,\tau,\mu)\to 
\TT=(\widetilde{T},{}^*\Cc,\Sigma,\widetilde{Q},\widetilde{\tau},\widetilde{\varphi}).$$
It is an isomorphism if $\Cc$ is locally projective as a
left $A$-module.
\end{proposition}

\begin{proof}
We have the inclusion
$$T=\End^\Cc(\Sigma)\subset \widetilde{T}=\End_{{}^*\Cc}(\Sigma)$$
We also have a map
$$\alpha:\ Q={}^\Cc\Hom(\Cc,\Sigma^*)\to \widetilde{Q}=
\Hom_{{}^*\Cc}(\Sigma,{}^*\Cc).$$
For $q:\ \Cc\to \Sigma^*$, we let
$$\alpha(q)={}^*q:\ {}^*(\Sigma^*)\cong \Sigma\to {}^*\Cc;$$
the fact that $\Sigma$ is isomorphic to its double dual follows from the
fact that $\Sigma$ is finitely generated and projective as a right $A$-module.
Let us show that left $\Cc$-colinearity of $q$ implies right ${}^*\Cc$-linearity
of ${}^*q$. First observe that ${}^*q(u)(c)=q(c)(u)$. For all $f\in {}^*\Cc$,
$u\in \Sigma$ and $c\in \Cc$, we have
\begin{eqnarray*}
&&\hspace*{-2cm}
{}^*q(u\cdot f)(c)={}^*q(u_{[0]}f(u_{[1]}))(c)=
q(c)(u_{[0]}f(u_{[1]}))\\
&=& q(c)(u_{[0]})f(u_{[1]})=f(q(c)(u_{[0]})u_{[1]})\\
{\rm \equref{4.6.1}}~~~~
&=& f(c_{(1)}(q(c_{(2)})(u)))= f(c_{(1)}({}^*q(u)(c_{(2)})))\\
&=& ({}^*q(u)\# f)(c)
\end{eqnarray*}
Let us show that this defines a morphism of Morita contexts, i.e.
$$\mu=\widetilde{\mu}\circ (\alpha\ot \Sigma)~~{\rm and}~~
\tau=\widetilde{\tau}\circ (\Sigma\ot\alpha)$$
Indeed,
$$\widetilde{\mu}({}^*q\ot u)(c)={}^*q(u)(c)=q(c)(u)=\mu(q\ot u)(c)$$
and
$$\widetilde{\tau}(u\ot {}^*q)(v)=u_{[0]}({}^*q(v)(u_{[1]}))=
u_{[0]}(q(u_{[1]})(v))=\tau(u\ot q)(v).$$
Now assume that $\Cc\in {}_A\Mm$ is locally projective. Recall that this means
that, for every finite $D\subset \Cc$, there exists $\sum_i c_i^*\ot c_i\in
{}^*\Cc\ot_A\Cc$ such that
$$d=\sum_i c_i^*(d)c_i,$$
for all $d\in D$. We first show that $\widetilde{T}\subset T$. Take $f\in \widetilde{T}$,
and fix $u\in \Sigma$. Then write
\begin{equation}\eqlabel{4.8.1}
\rho(u)=u_{[0]}\ot u_{[1]}=\sum_{j=1}^n u_j\ot d_j,~~
\rho(f(u))=f(u)_{[0]}\ot f(u)_{[1]}=\sum_{k=1}^m v_k\ot e_k,
\end{equation}
and consider the finite set $D=\{d_1,\cdots,d_n,e_1,\cdots,e_m\}\subset \Cc$.
Taking $\sum_i c_i^*\ot c_i\in{}^*\Cc\ot_A\Cc$, as above, we can compute
\begin{eqnarray*}
&&\hspace*{-2cm}
f(u_{[0]})\ot_A u_{[1]}= \sum_i f(u_{[0]})\ot_A c_i^*(u_{[1]})c_i\\
&=& \sum_i f(u_{[0]}c_i^*(u_{[1]}))\ot_A c_i=\sum_i f(u\cdot c_i^*)\ot_A c_i\\
&=& \sum_i f(u)\cdot c_i^*\ot_A c_i=\sum_i f(u)_{[0]} c_i^*(f(u)_{[1]})\ot_A c_i\\ 
&=& \sum_i f(u)_{[0]} \ot_Ac_i^*(f(u)_{[1]})c_i=f(u)_{[0]}\ot_Af(u)_{[1]}
\end{eqnarray*}
proving that $f$ is right $\Cc$-colinear, as needed.\\
Now take $\widetilde{q}\in \widetilde{Q}$, and define $q=\beta(\widetilde{q})$
by
$$q(c)(u)=\widetilde{q}(u)(c)$$
We will show, using \leref{4.6}, that $q$ is right $\Cc$-colinear. We know that
$\widetilde{q}$ is right ${}^*\Cc$-linear, hence
$$\widetilde{q}(u\cdot f)(c)=f(c_{(1)}(\widetilde{q}(u)(c_{(2)}))),$$
and
\begin{equation}\eqlabel{4.8.2}
q(c)(u\cdot f)=f(c_{(1)}(q(c_{(2)})(u))),
\end{equation}
for all $c\in \Cc$, $f\in {}^*\Cc$ and $u\in \Sigma$. Fix $u\in \Sigma$, let
$D= \{d_1,\cdots,d_n\}\subset \Cc$ as in \equref{4.8.1}, and take the corresponding
$\sum_i c_i^*\ot c_i\in{}^*\Cc\ot_A\Cc$. We then compute
\begin{eqnarray*}
&&\hspace*{-2cm}
(q(c)(u_{[0]}))u_{[1]}=\sum_i q(c)(u_{[0]})c_i^*(u_{[1]})c_i\\
&=& \sum_i q(c)(u_{[0]}c_i^*(u_{[1]}))c_i=\sum_i q(c)(u\cdot c_i^*)c_i\\
{\rm \equref{4.8.2}}~~~~&=&\sum_i c_i^*(c_{(1)}(q(c_{(2)})(u)))c_i=
c_{(1)}(q(c_{(2)})(u))
\end{eqnarray*}
This proves that $q$ satisfies \equref{4.6.1}, hence $q\in Q$. We have a well-defined
map $\beta:\ \widetilde{Q}\to Q$, which is clearly the inverse of $\alpha$.
\end{proof}

Now let $A$ and $B$ be rings, and $\Sigma\in {}_B\Mm_A$. We will compare
the Morita context $\DD$ associated to the comatrix coring $\Dd=\Sigma^*\ot_B\Sigma$
to the Morita context $\SS$ associated to $\Sigma\in {}_B\Mm$. Recall that this
Morita context is
$$\SS=(B,S={}_B\End(\Sigma)^{\rm op},\Sigma, {}^*\Sigma={}_B\Hom(\Sigma,B),
\varphi,\psi)$$
with
$$\varphi:\ \Sigma\ot_S{}^*\Sigma\to B,~~\varphi(u\ot s)=s(u)$$
and
$$\psi:\ {}^*\Sigma\ot_B\Sigma\to S,~~\psi(\gamma\ot u)(v)=\gamma(u)v.$$

\begin{proposition}\prlabel{4.9}
With notation as above, we have a morphism of Morita contexts
$$\SS=(B,S,\Sigma, {}^*\Sigma,\varphi,\psi)\to
\DD=(T=\End^{\Dd}(\Sigma),{}^*\Dd,\Sigma,Q,\tau,\mu).$$
It is an isomorphism if $\Sigma\in {}_B\Mm$ is totally faithful.
\end{proposition}

\begin{proof}
If $\Sigma\in {}_B\Mm$ is totally faithful, then the map
$$\eta_N:\ N\to \Hom^\Dd(\Sigma,N\ot_B\Sigma),~~\eta(n)(u)=n\ot_B u$$
is an isomorphism, for every $N\in {}_B\Mm$. In particular,
$\eta_B:\ B\to T=\End^{\Dd}(\Sigma)$ is then an isomorphism. Since
$\Sigma\in \Mm_A$ is finitely generated projective, we also have an isomorphism
$${}^*\Dd={}_A\Hom(\Sigma^*\ot_B \Sigma,A)\cong S={}_B\End(\Sigma)$$
We will next construct a map
$$\lambda:\ {}^*\Sigma={}_B\Hom(\Sigma,B)\to Q={}^\Dd\Hom(\Dd,\Sigma^*).$$
A left $A$-linear map $\varphi:\ \Dd\to \Sigma^*$ belongs to $Q$
(i.e. is left $\Dd$-colinear) if and only if
\begin{equation}\eqlabel{4.9.1}
\sum_i \varphi(f\ot_B u)\ot_B e_i\ot_A f_i=
\sum_i f\ot_B e_i\ot_A \varphi(f_i\ot_B u).
\end{equation}
Take $\gamma\in {}^*\Sigma={}_B\Hom(\Sigma,B)$, and define
$\lambda(\gamma)=\varphi$ by
$$\varphi(f\ot_B u)=f\gamma(u).$$
then $\lambda(\gamma)\in Q$ since
$$f\gamma(u)\ot_Be_i\ot_A f_i=f\ot_B\gamma(u)e_i\ot_A f_i=
f\ot_Be_i\ot_A f_i\gamma(u).$$
If $\Sigma\in {}_B\Mm$ is totally faithful, then the inverse
$\ol{\lambda}$ of $\lambda$ is 
$$\ol{\lambda}:\ {}^\Dd\Hom(\Dd,\Sigma^*)\to {}_B\Hom(\Sigma, B)\cong
{}_B\Hom(\Sigma,\End^\Dd(\Sigma)),$$
given by $\ol{\lambda}(\varphi)=\beta$, with
\begin{equation}\eqlabel{4.9.2}
\beta(u)(v)=\sum_i e_i(\varphi(f_i\ot_B u)(v)).
\end{equation}
We prove that $\beta(u)\in \End^\Dd(\Sigma)$ it suffices to show that
\begin{equation}\eqlabel{4.9.3}
\sum_j \beta(u)(e_j)\ot_A f_j\ot_B v=\sum_j e_j\ot_A f_j \ot_B \beta(u)(v)
\end{equation}
or $A=B$, where
\begin{eqnarray*}
A&=& \sum_{i,j} e_i(\varphi(f_i\ot_B u)(e_j))\ot_A f_j\ot_B v\\
B&=& \sum_{i,j} e_j\ot_A f_j\ot_B e_i(\varphi(f_i\ot_B u)(v))
\end{eqnarray*}
It follows from \equref{4.9.1} that
$$\sum_{i,j} e_i\ot_A \varphi(f_i\ot_B u)\ot_B e_j\ot_A f_j\ot_B v=
\sum_{i,j}e_i\ot_A f_i\ot_B e_j\ot_A \varphi(f_j\ot_B u)\ot_B v,$$
and, after we let the second tensor factor act on the third one,
$$A= \sum_j e_j\ot_A \varphi(f_j\ot_B u)\ot_B v$$
Using \equref{4.9.1}, we also obtain that
$$\sum_{i,j} e_j\ot_A f_j\ot_B e_i\ot_A \varphi(f_i\ot_B u)\ot_B v=
\sum_{i,j} e_j\ot_A \varphi(f_j\ot_B u)\ot_B e_i\ot_A f_i\ot_B v;$$
letting the fourth tensor factor act on the fifth, we find
$$B= \sum_j e_j\ot_A \varphi(f_j\ot_B u)\ot_B v,$$
and \equref{4.9.3} follows.\\
Let us now check that $\lambda$ and $\ol{\lambda}$ are inverses, at
least if we identify $B$ and $T$. Take $\lambda\in {}_B\Hom(\Sigma,B)$,
and $(\ol{\lambda}\circ\lambda)(\gamma)=\beta:\ \Sigma\to \End^\Dd(\Sigma)$.
Then
$$\beta(u)(v)=\sum_i e_i(f_i\gamma(u))(v)=\sum_i e_i(f_i(\gamma(u)(v)=
\gamma(u)v,$$
as needed. Now take $\varphi\in {}^\Dd\Hom(\Dd,\Sigma^*)$, and put
$\beta=\ol{\lambda}(\varphi)$, $\psi=\lambda(\beta)$. Then
\begin{eqnarray*}
&&\hspace*{-2cm}
\psi(f\ot_B u)(v)=f(\beta(u)(v))=f(\sum_i e_i(\varphi(f_i\ot_B u)(v)))\\
&=& \sum_i f(e_i)(\varphi(f_i\ot_B u)(v))=\varphi(f\ot_B u)(v)
\end{eqnarray*}
To show that we really have a morphism of Morita contexts, we first have
to show that the diagram
$$\begin{diagram}
\Sigma\ot_S{}^*\Sigma&\rTo^{\varphi}&B\\
\dTo^{\Sigma\ot_S\lambda}&&\dTo_{\eta_B}\\
\Sigma\ot_{{}^*\Dd}Q&\rTo^{\tau}&T
\end{diagram}$$
commutes. Indeed,
\begin{eqnarray*}
&&\hspace*{-2cm}
\Bigl((\tau\circ (\Sigma\ot_S\lambda))(u\ot \gamma)\Bigr)(v)=
\tau(u\ot\lambda(\gamma))(v)\\
&=& u_{[0]}(\lambda(\gamma)(u_{[1]}))(v)
= \sum_i e_i(\lambda(\gamma)(f_i\ot_B u))(v)\\
&=& \sum_i e_if_i(\gamma(u)(v))
= \gamma(u)(v)= ((\eta_B\circ\varphi)(u\ot\gamma))(v)
\end{eqnarray*}
Finally, we need commutativity of the diagram
$$\begin{diagram}
{}^*\Sigma\ot_B\Sigma&\rTo^{\psi}&S\\
\dTo^{\lambda\ot\Sigma}&&\dTo_{\alpha}\\
Q\ot_T\Sigma&\rTo^{\mu}&{}^*\Dd
\end{diagram}$$
This is also straightforward:
\begin{eqnarray*}
&&\hspace*{-2cm}
\Bigl(\bigl(\tau\circ (\Sigma\ot_S\lambda)\Bigl(u\ot \gamma\Bigr)(v)=
(\tau(u\ot \lambda(\gamma))(v)\\
&=& u_{[0]}(\lambda(\gamma)(u_{[1]}))(v)
= \sum_i e_i (\lambda(\gamma)(f_i\ot u))(v)\\
&=& \sum_i e_if_i(\gamma(u)v)=\gamma(u)v\
= \Bigl((\eta_B\circ \varphi)(u\ot\gamma)\Bigr)(v)
\end{eqnarray*}
\end{proof}

\begin{proposition}\prlabel{4.10}
Consider the Morita context $\CC=(T,{}^*\Cc,{}^*M,Q,\tau,\mu)$
from \thref{4.4}, and assume that $M$ is finitely generated and projective as a left 
$A$-module. The following statements are equivalent:
\begin{enumerate}
\item $\tau$ is surjective (hence bijective);
\item for every $N\in \Mm^\Cc$, the map
$$\omega_N:\ N\ot_{{}^*\Cc} Q\to \Hom^\Cc({}^*M,N),~~
\omega_N(n\ot q)(u)=n\mu(q\ot u)$$
is surjective;
\item the natural transformation
$$\omega:\ \bullet\ot_{{}^*\Cc} Q\to \Hom^\Cc({}^*M,\bullet)$$
given by
$$\omega_N:\ N\ot_{{}^*\Cc} Q\to \Hom^\Cc({}^*M,N),~~
\omega_N(n\ot q)(u)=n\mu(q\ot u),$$
for every $N\in \Mm^\Cc$, is an isomorphism;
\item the natural transformation
$$\ol{\omega}:\ \tilde{G}=\bullet\ot_{{}^*\Cc} Q\to G= \bullet\otimes^\Cc{M}$$
given by
$$\ol{\omega}_N:\ N\ot_{{}^*\Cc} Q\to N\ot^\Cc M,~~
\ol{\omega}_N(n\ot q)=\sum_i n\mu(q\ot f_i)\ot e_i,$$
for every $N\in \Mm^\Cc$, is an isomorphism.
\end{enumerate}
In this case, $M$ is finitely generated and projective as a left $T$-module.
\end{proposition}

\begin{proof}
$1)\Rightarrow 3)$. If $\tau$ is surjective, then 
${}^*M$ is finitely generated and projective as a right $T$-module
(\cite[Theorem I.3.4]{Ba}), so $M$ is finitely generated and projective
as a left $T$-module.\\
Take $\Sigma={}^*M$,
and let $\sum_i e_i\ot_A f_i$ be a finite dual basis of $\Sigma\in\Mm_A$, as
before.
Choose $u_j\in {}^*M$ and $q_j\in Q$ such that
$\tau(\sum_j u_j\ot q_j)$ is the identity map on $M$. Then we define
$$\psi_N:\ \Hom^\Cc({}^*M,N)\to N\ot_{{}^*\Cc} Q,~~\psi_N(\varphi)=
\sum_j \varphi(u_j)\ot q_j.$$
Then $\psi_N$ and $\omega_N$ are inverses:
\begin{eqnarray*}
&&\hspace*{-2cm} \psi_N(\omega_N(n\ot_{{}^*\Cc} q))=
\sum_j n\mu(q\ot u_j)\ot_{{}^*\Cc} q_j\\
&=& \sum_j n\ot_{{}^*\Cc}\mu(q\ot u_j) q_j=
 \sum_j n\ot_{{}^*\Cc}q\tau(u_j\ot q_j)= n\ot_{{}^*\Cc}q
\end{eqnarray*}
and
\begin{eqnarray*}
&&\hspace*{-2cm} 
\omega_N(\psi_N(\varphi))(u)= \omega_N(\sum_j \varphi(u_j)\ot q_j)(u)\\
&=& \sum_j \varphi(u_j)\mu(q_j\ot u)
= \sum_j \varphi(u_j)_{[0]}(q_j( \varphi(u_j)_{[1]})(u))\\
&=& \sum_j \varphi(u_{j[0]})(q_j( u_{j[1]})(u))
= \sum_i \varphi(e_i)(f_i(u))=\varphi(u).
\end{eqnarray*}
$3)\Rightarrow 2)$ is trivial.\\
$2)\Rightarrow 1)$: take $N=M$.
$3)\Leftrightarrow 4)$ follows from \prref{1.3}.
\end{proof}

From now on, we restrict attention to the Morita context $\Cc=(T=\End^\Cc(\Sigma), {}^*\Cc,\Sigma,Q={}^\Cc\Hom(\Cc,\Sigma^*),
\tau,\mu)$ from \equref{4.5.1}, with $\Sigma\in \Mm_A$
finitely generated projective. We study the image of the map $\mu$. Assume that
$\Cc\in {}_A\Mm$ is locally projective, and recall from \cite{CVW2} that
$f\in {}^*\Cc$ is called rational if there exist a finite number $f_i\in \Cc^*$
and $c_i\in \Cc$ such that
\begin{equation}\eqlabel{4.10a.1}
f\# g= \sum_f f_i g(c_i)
\end{equation}
for all $g\in {}^*\Cc$. Then
$$({}^*\Cc)^{\rm rat}=\{f\in {}^*\Cc~|~f~{\rm~is~rational}\}$$
is a right $\Cc$-comodule. 

\begin{lemma}\lelabel{4.10a}
Let $\Sigma\in \Mm^\Cc_{\rm fgp}$, where $\Cc$ is locally projective as a left
$A$-module, and consider the $\mu$ from the Morita context \equref{4.5.1}.
Then
$$\im\mu\subset ({}^*\Cc)^{\rm rat}.$$
\end{lemma}

\begin{proof}
Take $\mu(q\ot u)\in \im\mu$.
For all $f\in {}^*\Cc$ and $c\in \Cc$, we have
\begin{eqnarray*}
&&\hspace*{-2cm}
(\mu(q\ot u)\# f)(c)=f(c_{(1)}\mu(q\ot u)(c_{(2)}))\\
&=& f(c_{(1)}q(c_{(2)})(u))=f(q(c)(u_{[0]})u_{[1]})\\
&=& q(c)(u_{[0]})f(u_{[1]})=\mu(q\ot u_{[0]})(c)f(u_{[1]}),
\end{eqnarray*}
and the rationality of $\mu(q\ot u)$ follows after we take $f_i=\mu(q\ot u_{[0]})$
and $c_i=u_{[1]}$.
\end{proof}

\begin{corollary}\colabel{4.10b}
If $\mu$ is surjective, then $\Cc$ is finitely generated and projective as a left $A$-module.
\end{corollary}

\begin{proof}
If $\mu$ is surjective, then it follows from \leref{4.10a} that every $f\in {}^*\Cc$
is rational, and then it follows from \cite[Cor. 4.2]{CVW2} that $\Cc\in {}_A\Mm$
is finitely generated projective.
\end{proof}

We now consider the situation where $\Cc$ is finitely generated and projective
as a left $A$-module. Then the categories $\Mm^\Cc$ and $\Mm_{{}^*\Cc}$ are
isomorphic. The functor
$$F=\bullet\ot_B\Sigma:\ \Mm_B\to \Mm^\Cc\cong \Mm_{{}^*\Cc}$$
has a right adjoint  $G=\Hom^\Cc(\Sigma,\bullet)$. If the map $\tau$ in the Morita
context $\CC$ is surjective, and $B=T$,
then $\tilde{G}=\bullet\ot_{{}^*\Cc}Q$ is also a right
adjoint of $F$, hence $G\cong \tilde{G}$, by Kan's Theorem. If we construct the
isomorphism $G(M)\cong \tilde{G}(M)$, following for example \cite[Prop. 9]{CMZ},
then we recover the isomorphism from \prref{4.10}.\\
We are now able to state and prove the main result of this Section. It generalizes
\cite[Theorem 4.7]{Caenepeel03}.

\begin{theorem}\thlabel{4.11}
Let $A$ and $B$ be rings, and $\Cc$ an $A$-coring, which is finitely generated and
projective as a left $A$-module. Let $\Sigma\in {}_B\Mm^{\Cc}_{\rm fgp}$. Also
write $T=\End^\Cc(\Sigma)$. 
Then the following assertions are equivalent:
\begin{enumerate}
\item
\begin{itemize}
\item ${\rm can}:\ \Dd=\Sigma^*\ot_B \Sigma\to \Cc$ is an isomorphism;
\item $\Sigma\in {}_B\Mm$ is faithfully flat.
\end{itemize}
\item
\begin{itemize}
\item ${}^*{\rm can}:\ {}^*\Cc\to {}_B\End(\Sigma)^{\rm op}$ is an isomorphism;
\item $\Sigma\in {}_B\Mm$ is progenerator.
\end{itemize}
\item
\begin{itemize}
\item $l:\ B\to T$ is an isomorphism;
\item the Morita context $\CC=(T,{}^*\Cc,\Sigma,Q,\tau,\mu)$ from \equref{4.5.1}
is strict. 
\end{itemize}
\item
\begin{itemize}
\item $(F,G)$ is a pair of inverse equivalences between the categories
$\Mm_B$ and $\Mm^\Cc$.
\end{itemize}
\end{enumerate}
\end{theorem}

\begin{proof}
$\ul{1)\Rightarrow 4)}$. From the faithfully flat descent \thref{2.7}, $(K,R)$ is
a pair of equivalences; the fact that ${\rm can}$ is an isomorphism then implies
that $(F,G)$ is an isomorphism, by \prref{3.2}.\\
$\ul{4)\Rightarrow 2)}$. $F=\bullet_B\Sigma$ is an equivalence between the
module categories $\Mm_B$ and $\Mm_{{}^*\Cc}$, hence $\Sigma$ is a left
$B$-progenerator. It follows from \prref{3.3} that ${\rm can}$ is an isomorphism,
and then the dual map ${}^*{\rm can}$ is also an isomorphism.\\
$\ul{2)\Rightarrow 1)}$. It follows from \prref{4.1} that ${\rm can}$ is an isomorphism.\\
$\ul{4)\Rightarrow 3)}$. It follows from \prref{3.3} that $l$ is an isomorphism. 
since 4) implies 2), we know that $\Sigma\in {}_B\Mm$ is a progenerator. Then the
associated Morita context $\SS$ is strict. The Morita context $\DD$ is then also strict,
since it is isomorphic to it (see \prref{4.9}). Now 4) implies 1), so can is an isomorphism,
and $\Cc$ is isomorphic to $\Dd$ as a coring, hence $\CC\cong \DD$ is also strict.\\
$\ul{3)\Rightarrow 4)}$. If $\CC$ is strict, then $F$ is an equivalence of categories.
\end{proof}

We now look at the situation where $\Cc$ is locally projective as a left $A$-module.
If $R$ is a ring with local units, then we denote by $\Mm_R$ the category of
right unital $R$-modules, these are right $R$-modules for which the canonical
map $M\ot_RR\to R$ is an isomorphism.

\begin{lemma}\lelabel{4.12}
Let $\Cc$ be an $A$-coring which is locally projective as a left $A$-module.
The rational dual $({}^*\Cc)^{\rm rat}$ has local units if and only if
$({}^*\Cc)^{\rm rat}$ is dense in ${}^*\Cc$ with respect to the finite topology.
In this situation, we have the following properties.
\begin{enumerate}
\item For every $M\in \Mm^\Cc$, the map
$$\Omega_M:\ M\ot_{{}^*\Cc}({}^*\Cc)^{\rm rat}\to M,~~
\Omega_M(m\ot_{{}^*\Cc}f)=m\cdot f=m_{[0]}f(m_{[1]})$$
is an isomorphism.
\item The categories $\Mm_{({}^*\Cc)^{\rm rat}}$ and $\Mm^\Cc$ are isomorphic.
\end{enumerate}
\end{lemma}

\begin{proof}
For the first statement, we refer to \cite[Prop. 4.1]{CVW2}.\\
1) We define $\Psi_M:\ M\to M\otimes_{\*C}(\*C)^{\rm rat}$, $\Psi_M(m)=m\otimes e$, where $e$ (depending on $m$) is constructed as follows. Write $\rho(m)=\sum_jm_j\otimes_A c_j$. Then pick $e\in(\*C)^{\rm rat}$ such that $\varepsilon(c_j)=e(c_j)$, for all $j$. This means $e$ acts as a local unit on $m$: 
$$m\cdot e=m_je(c_j)=m_j\varepsilon(c_j)=m_{[0]}\varepsilon(m_{[1]})=m.$$

We first check that $\Psi_M$ is well-defined. Take another $e'\in(\*C)^{\rm rat}$ satisfying $\varepsilon(c_j)=e'(c_j)$. We have to show $m\otimes_{\*C} e=m\otimes_{\*C} e'$. To this end, choose any common local unit $e''\in (\*C)^{\rm rat}$ for $e$ and $e'$, i.e. $e=ee''$ and $e'=e'e''$, then we compute
\begin{eqnarray*}
&&\hspace*{-2cm}
m\otimes_{\*C}e=m\otimes_{\*C}ee''=m\cdot e\otimes_{\*C}e''=
m \otimes_{\*C}e''\\
&=&m\cdot e' \otimes_{\*C}e''=m\otimes_{\*C}e'e''=m\otimes_{\*C}e'.
\end{eqnarray*}
$\Omega_M$ is a left inverse of $\Psi_M$, since
$$\Omega_M(\Psi_M(m))=m\cdot e=m.$$
To show that $\Omega_M$ is a right inverse of $\Psi_M$,
take $m\otimes f\in M\otimes_{\*C}(\*C)^{\rm rat}$. Write $\Psi(m\cdot f)=m\cdot f\otimes e$, and pick a common local unit $e'\in (\*C)^{\rm rat}$ for $f$ and $e$.
\begin{eqnarray*}
&&\hspace*{-2cm}
\Psi_M(\Omega_M(m\otimes_{\*C} f))=m\cdot f\otimes_{\*C} e
=m\cdot f\otimes_{\*C} ee'\\
&=&(m\cdot f)\cdot e\otimes_{\*C} e'=m\cdot f\otimes_{\*C} e'\\
&=&m\otimes_{\*C} fe'=m\otimes_{\*C} f.
\end{eqnarray*}
2) Starting with $M\in\Mm^\Cc$, we find $M\in\Mm_{\*C}$ and as in 1) one shows that the restricted action of $(\*C)^{\rm rat}$ on $M$ is unital, so $M\in\Mm_{(\*C)^{\rm rat}}$. Conversely, if $M\in\Mm_{(\*C)^{\rm rat}}$, then for every $m\in M$ we can find elements $m_i\in M$ and $g_i\in(\*C)^{\rm rat}$ such that $m=m_i\cdot g_i$. For all $f\in\*C$ we then compute
$$m\cdot f=(m_i\cdot g_i)\cdot f=m_i\cdot(g_i\#f)=m\cdot g_{i[0]}f(g_{i[1]}).$$
This means that $M$ is a rational $\*C$-module, hence $M\in\Mm^{\Cc}$.
\end{proof}

We now present a generalization of \leref{4.10a}.

\begin{corollary}\colabel{4.14}
Consider an $A$-coring $\Cc$ which is locally projective as left $A$-module and let $\mu$ be as in the Morita context from \thref{4.4}. Then for every $M\in\Mm_{\*C}$ we have a map
$$r_M:\ \tilde{F}\tilde{G}M=M\otimes_{\*C}Q\otimes_T\Sigma\rTo^{I_M\otimes\mu}M^{\rm rat},$$
which is an isomorphism if $\im\mu=(\*C)^{\rm rat}$ and $(\*C)^{\rm rat}$ has right local units.
\end{corollary}

\begin{proof}
First of all, $r_M$ is well defined: pick $m\otimes q\otimes u\in M\otimes_{\*C}Q\otimes_T\Sigma$, then $r_M(m\otimes q\otimes u)=m\cdot\mu(q\otimes u)$. Since $\im\mu\subset(\*C)^{\rm rat}$, we find
\begin{eqnarray*}
&&\hspace*{-2cm}
(m\cdot\mu(q\otimes u))\cdot f=m\cdot(\mu(q\otimes u))\cdot f)=m\cdot((\mu(q\otimes u)\# f)\\
&=&\sum_im\cdot((\mu(q\otimes u)_if(c_i))=
\sum_i(m\cdot(\mu(q\otimes u)_i)f(c_i),
\end{eqnarray*}
so we conclude that $m\cdot\mu(q\otimes u)\in M^{\rm rat}$.\\
If $(\*C)^{\rm rat}$ has local units, then, as 
$M^{\rm rat}\in\Mm^\Cc$, $M^{\rm rat}\cong M\otimes_{\*C} (\*C)^{\rm rat}$, by 
\leref{4.12}. If, in addition, $\im\mu=(\*C)^{\rm rat}$, then this isomorphism is exactly $r_M$.
\end{proof}

\coref{4.14} provides an explicit way to construct the rational part of a $\*C$-module. Remark that $r_{\*C}=\mu$.

We have seen that the Morita context $\CC=(T,\*C,\Sigma,Q,\tau,\mu)$ can only be strict if $\Cc$ is finitely generated and projective as a left $A$-module, by the surjectivity of $\mu$. Consequently, in many cases, it is better to look to an other, restricted, Morita context. Since $\im\mu\subseteq (\*C)^{\rm rat}$, we can restrict our context without any consequenses on the connecting maps or modules to $\CC'= (T,(\*C)^{\rm rat},\Sigma,Q,\tau,\mu)$. 

If $(\*C)^{\rm rat}$ satisfies the conditions of \leref{4.12}, then we have a Morita context connecting the ring with unit $T$ and the ring with local units $(\*C)^{\rm rat}$. This has the following implications (for details see \cite[Prop 2.12]{CVW2} and \cite{AM}):
\begin{enumerate}
\item the bijectivity of $\mu$ and $\tau$ follows from their surjectivity;
\item if $\tau$ is surjective, then $\Sigma_{(\*C)^{\rm rat}}$,
$\Sigma_{\*C}$, ${_{(\*C)^{\rm rat}}Q}$ and
${_{\*C}Q}$ are finitely generated and projective (using the Morita contexts
$\CC'$ and $\CC$);
\item if $\im\mu=(\*C)^{\rm rat}$, then ${}_T\Sigma$ and $Q_T$ are locally projective.
\end{enumerate}

\begin{theorem}\thlabel{4.15}
Let $\Cc$ be an $A$ coring which is locally projective as left $A$-module. Suppose $(\*C)^{\rm rat}$ is  dense in the finite topology on $\*C$. Take $\Sigma\in\Mm^\Cc_{\rm fgp}$ and let $\CC'=(T,(\*C)^{\rm rat}),\Sigma,Q,\mu,\tau)$ be the restricted morita context. If $\ell:\ B\to T$ is an isomorphism and $\tau$ is surjective, then the following statements are equivalent.
\begin{enumerate}
\item ${\rm can} :\ \Dd=\Sigma^*\otimes_B\Sigma\to\Cc$ is an isomorphism and ${_B\Sigma}$ is faithfully flat;
\item ${\rm can} :\ \Dd=\Sigma^*\otimes_B\Sigma\to\Cc$ is an isomorphism and ${_B\Sigma}$ is flat;
\item $\Sigma$ is a generator in $\Mm^\Cc$;
\item $\Sigma$ is a projective generator in $\Mm^\Cc$;
\item $\Sigma$ is a progenerator in $\Mm_{(\*C)^{\rm rat}}$;
\item $\mu$ is surjective (onto $(\*C)^{\rm rat}$);
\item $\CC'$ is a strict morita context;
\item $(F,G)$ is a pair of inverse equivalences between $\Mm_B$ and $\Mm^\Cc$;
\item for all $ N \in \Mm^\Cc$, the counit of the adjunction $\xi_N:\ \Hom^\Cc(\Sigma,N)\otimes_B\Sigma\to N$ is an isomorphism.
\end{enumerate}
\end{theorem}

\begin{proof}
$\ul{(1)\Leftrightarrow(4)\Leftrightarrow(8)}$ follow from
 \thref{3.9} and the fact that local projectivity implies flatness.\\
$\ul{(2)\Leftrightarrow(3)\Leftrightarrow(9)}$ follow in the same way
from \thref{3.8}.\\
$\ul{(6)\Leftrightarrow(7)}$ follows from Morita theory.\\
$\ul{(7)\Rightarrow (8)}$. Since $(\*C)^{\rm rat}$ is dense, $\Mm_{(\*C)^{\rm rat}}\cong\Mm^\Cc$ by \leref{4.12}. The strictness of the Morita context $\CC'$ implies 
that the categories $\Mm_T\cong\Mm_B$ and $\Mm_{(\*C)^{\rm rat}}\cong\Mm^\Cc$ 
are equivalent via $F$, see \prref{4.10}, (4).\\
$\ul{(8)\Rightarrow (9)}$ is trivial.\\
$\ul{(9)\Rightarrow (7)}$. Since $\Mm^\Cc$ is a full subcategory of $\Mm_{\*C}$, we
have that 
$$\Hom^\Cc(\Sigma,(\*C)^{\rm rat})\cong\Hom_{\*C}(\Sigma,(\*C)^{\rm rat})=\Hom_{\*C}(\Sigma,\*C)\cong Q.$$
Indeed, for all $\varphi\in\Hom_{\*C}(\Sigma,\*C)$, $f\in\*C$ and $u\in \Sigma$,
we have that
$$\varphi(u)\cdot f = \varphi(u\cdot f) = \varphi(u_{[0]}f(u_{[1]})) = \varphi(u_{[0]})f(u_{[1]}),$$
so we conclude that $\varphi(u)\in(\*C)^{\rm rat}$.\\
Now take $N=(\*C)^{\rm rat}$ in the counit of the adjunction; we then find that $\zeta_{(\*C)^{\rm rat}}=\mu$ is an isomorphism, as
$$\zeta_{(\*C)^{\rm rat}}:
\Hom_{\*C}(\Sigma,(\*C)^{\rm rat})\otimes_B\Sigma 
\cong Q\otimes_B \Sigma\to N =(\*C)^{\rm rat}.$$
$\ul{(4)\Rightarrow (5)}.$ We know $\Sigma$ is a generator in $\Mm^\Cc\cong\Mm_{(\*C)^{\rm rat}}$. Since (4) is equivalent to (7), we know that $\CC'$ is strict. From Morita theory it then follows that $\Sigma\in\Mm_{(\*C)^{\rm rat}}$
is finitely generated projective.\\
$\ul{(5)\Rightarrow (4)}$ is trivial.
\end{proof}

\section{Coseparable corings and an affineness Theorem}\selabel{5}
Let $A$ be a ring, $\Cc$ an $A$-coring, $\Sigma\in \Mm^\Cc_{\rm fgp}$ and $T=\End^{\Cc}(\Sigma)$. Then $\Sigma \in {}_T\Mm^\Cc_{\rm fgp}$ and we can consider the adjoint pairs of functors $(F,G)$ and $(F',G')$ introduced in \seref{1}. We also consider the comatrix coring $\Dd=\Sigma^*\ot_T\Sigma$. As we have seen in \seref{3}, $\Cc$ is Galois if the canonical map is bijective. In this section we will discuss when surjectivity of the canonical map is a sufficient condition for  $(F,G)$ and $(F',G')$ being a pair of inverse equivalences, and,
a fortiori, $(\Cc,\Sigma)$ being Galois. Properties of this type have been studied in
special situations in \cite{Militaru,SchSch,Schneider90}.

Recall \cite{CMZ} that we have two pairs of adjoint functors
$(H,Z)$ and $(H',Z')$, 
$$H :\ \Mm^\Cc \to \Mm_A~~;~~
Z =\bullet\ot_A \Cc:\ \Mm_A \to \Mm^\Cc$$
$$H' :\ {}^\Cc\Mm \to {}_A\Mm~~;~~
Z' =\Cc \ot_A\bullet:\ {}_A\Mm \to {}^\Cc\Mm$$
$H$ and $H'$ are the functors forgetting the $\Cc$-coaction. We have a bijective
correspondence between
$$V = \dul{\rm Nat}(ZH, 1_{\Mm^\Cc}),~~
\tilde{V} = \dul{\rm Nat}(Z'H', 1_{{}^\Cc\Mm})$$
and
$$V_2 = \{\theta \in {}_A\Hom_A(\Cc\ot_A \Cc,A) | c_{(1)}\theta(c_{(2)}\ot_A d)=\theta(c\ot_A d_{(1)})d_{(2)} \}.$$
We describe the correspondence between $V$ and $V_2$.If $\alpha \in V$ 
be a natural transformation, then $\theta = \alpha_{\Cc}\in V_2$.
Conversely, given $\theta \in V_2$, we define a natural transformation $\alpha$ by 
$$\alpha_N :\ N\ot_A \Cc \to N :\ n\ot_A c \to n_{(0)}\theta(n_{(1)}\ot_A c),$$
for all $N\in \Mm^\Cc$.

\begin{proposition}\prlabel{5.1}
Take $\theta \in V_2$, and let $\alpha \in V$ and $\beta \in \tilde{V}$ be the corresponding natural transformations. Then the following stements are equivalent
\begin{enumerate}
\item $\alpha_{\Sigma} \circ \rho^r = {\Sigma}$;
\item $\beta_{\Sigma^*} \circ\rho^l = {\Sigma^*}$;
\item $u_{(0)}\theta(u_{(1)}\ot_A u_{(2)}) = u$, for all $u \in \Sigma$;
\item $\theta(g_{(-2)}\ot_A g_{(-1)})g_{(0)} = g$, for all $g \in \Sigma^*$.
\end{enumerate}
\end{proposition}

\begin{proof}
We prove $4)\Rightarrow 3)$. The proof of the other applications is straightforward, and is left to the reader.
\begin{eqnarray*}
u &=&\sum_i e_if_i(u)=
\sum_i e_i\theta(f_{i[-2]} \ot_A f_{i[-1]})f_{i[0]}(u)\\
&=&\sum_{i,j} e_i\theta(f_i(e_{j[0]})e_{j[1]} \ot_A e_{j[2]})f_j(u)=
\sum_{i,j} e_if_i(e_{j[0]})\theta(e_{j[1]} \ot_A e_{j[2]})f_j(u)\\
&=&\sum_j e_{j[0]}\theta(e_{j[1]} \ot_A e_{j[2]})f_j(u)=
u_{[0]}\theta(u_{[1]} \ot_A u_{[2]})
\end{eqnarray*}
\end{proof}

$\theta$ is called $\Sigma$-normalized if the conditions
 of \prref{5.1} are satisfied.

\begin{lemma}\lelabel{5.t}
Assume $\theta \in V_2$ is $\Sigma$-normalized. Then there exists a surjective projection $t:\ \Sigma \ot_A \Sigma^* \to T = \Sigma \ot^{\Cc} \Sigma^*$ in ${}_T\Mm_T$.
\end{lemma}

\begin{proof}
We define
$$t(u \ot_A g) = u_{[0]}\theta(u_{[1]}\ot_A g_{[-1]})\ot_A g_{[0]}$$
From the fact that $\theta\in V_2$, it follows that
$$u_{[0]}\ot u_{[1]}\theta(u_{[2]}\ot_A g_{[-1]})\ot_A g_{[0]} = 
u_{[0]}\ot_A \theta(u_{[1]}\ot_A g_{[-2]})g_{[-1]}\ot_A g_{[0]}$$
which means precisely that 
$t(u\ot_A g) \in T=\Sigma \ot^{\Cc} \Sigma^*$.
Now take $u \ot_A g \in B = \Sigma \ot^{\Cc} \Sigma^*$. Then we have that
$$u_{[0]} \ot_A u_{[1]} \ot_A g  =  u \ot_A g_{[-1]} \ot_A g_{[0]}$$
hence
$$u_{[0]}\ot_A u_{[1]}\ot_A u_{[2]}\ot_A g = u_{[0]}\ot_A u_{[1]}\ot_A g_{[-1]}\ot_A g_{[0]}$$
and
$$u_{[0]}\theta(u_{[1]}\ot_A u_{[2]})\ot_A g = u_{[0]}\theta(u_{[1]}\ot_A g_{[-1]})\ot_A g_{[0]}$$
so
 $u \ot_A g = t(u \ot_A g)$.
A straightforward verification shows that $t$ is a morphism of $(T\hbox{-}T)$-bimodules.
\end{proof}

\begin{theorem}\thlabel{5.2}
Assume $\theta \in V_2$ is $\Sigma$-normalized. Then
$$\nu_N:\ N\to (N\ot_B\Sigma)\ot^\Cc\Sigma^*,~~
\nu_N(n)=\sum_i (n\ot_B e_i)\ot_A f_i,$$
is an isomorphism of right $B$-modules for all $N \in \Mm_B$. Hence $F$ is a fully faithful functor.
\end{theorem}

\begin{proof}
The inverse of $\nu_N$ is defined by
$$\theta_N :\ (N\ot_B \Sigma) \ot^{\Cc} \Sigma^* \to N :\ \sum_i (n_i \ot_B u_i) \ot_A g_i \to \sum_i n_i\cdot t(u_i\ot_A g_i)$$
Indeed, for all $n\in N$, we have that
$$
\theta_N \circ \nu_N (n) = \theta_N(\sum_i (n\ot_B e_i) \ot_A f_i) =
\sum_i n\cdot t(e_i\ot_A f_i) =
n\cdot e = n.$$
For all $\sum_j (n_j \ot_B u_j) \ot_A g_j\in (N\ot_B \Sigma) \ot^{\Cc} \Sigma^* $, we have
that
\begin{eqnarray*}
&&\hspace*{-2cm}
(\nu_N \circ \theta_N) (\sum_j (n_j \ot_B u_j ) \ot_A g_j ) = \nu_N(\sum_j n_j \cdot t(u_j \ot_A g_j )) \\
&=& \sum_{i,j} \left((n_j \ot_B e_i) \ot_A f_i \right)\cdot t(u_j \ot_A g_j ) \\
&=&
 \sum_{i,j} \left((n_j \ot_B e_i) \ot_A f_i \right)\cdot u_{j[0]}\theta(u_{j[1]}\ot_A g_{j[-1]})\ot_A g_{j[0]} \\
&=& \sum_{i,j} (n_j \ot_B e_i) \ot_A f_i(u_{j[0]})\theta(u_{j[1]}\ot_A g_{j[-1]})\ot_A g_{j[0]} \\
&=& \sum_j n_j \ot_B u_{j[0]}\theta(u_{j[1]}\ot_A g_{j[-1]})\ot_A g_{j[0]} \\
&=&\sum_j n_j \ot_B u_{j[0]}\theta(u_{j[1]}\ot_A u_{j[2]})\ot_A g_j  =
\sum_j (n_j \ot_B u_j) \ot_A g_j
\end{eqnarray*}
In the sixth equality, we used the fact that 
$\sum_j (n_j \ot_B u_j) \ot_A g_j\in (N\ot_B \Sigma) \ot^{\Cc} \Sigma^* $.
\end{proof}

\begin{lemma}\lelabel{5.3}
Assume that $\theta \in V_2$ is $\Sigma$-normalized. Then $\zeta_{F(P)}=\zeta_{P\ot_B \Sigma}$
is an isomorphism of right $\Cc$-comodules, for every $P \in \Mm_B$.
\end{lemma}

\begin{proof}
Since $(F,G)$ are an adjoint pair, we have that
$${F(P)} = \zeta_{F(P)} \circ F(\nu_P)$$
or
$${P\ot_B \Sigma} = \zeta_{P\ot_B \Sigma} \circ (\nu_P\ot_B {\Sigma})$$
We know from \leref{5.3} that $\nu_P$ is an isomorphism, so it follows that $\zeta_{P\ot_B \Sigma}$ is also an isomorphism.
\end{proof}

Recall that $\Cc$ is called a coseparable coring if the forgetful functor $H$ (and $H'$)
are separable. Recall from \cite{CMZ} that this is equivalent to the existence of
a natural transformation $\alpha\in V$ such that $\alpha\circ \eta$ is the identity
natural transformation, that is,
$$\alpha_{\Sigma}\circ\eta_{\Sigma} =\alpha_{\Sigma}\circ\rho^r_{\Sigma} = {\Sigma},$$
for all $\Sigma \in \Mm^{\Cc}$. Let $\theta\in V_2$ be the corresponding map. Then it
follows that $\theta$ is $\Sigma$-normalized, for every $\Sigma\in \Mm^\Cc$.

\begin{lemma}\lelabel{5.4}
Let $\Cc$ be a coseparable coring. Then there exists a $\theta \in V_2$ such that $\theta$ is $\Sigma$-normalized for every $\Cc$-comodule $\Sigma$.
\end{lemma}

\begin{proposition}\prlabel{5.5}
Assume that
\begin{enumerate}
\item $\Cc$ is projective as a right $A$-module;
\item $\Cc$ is a coseparable coring;
\item ${\rm can}$ is surjective.
\end{enumerate}
Then $(F,G)$ is a pair of inverse equivalences.
\end{proposition}

\begin{proof}
Taking into account \thref{5.2}, we only have to prove that $\zeta_M$ is an isomorphism, for all $M \in \Mm^{\Cc}$.
The map $\rho_M:\ M \to M \ot_A \Cc$ in $\Mm^{\Cc}$ has a left inverse $I_M \ot_A \varepsilon_{\Cc}$ in $\Mm_A$. The forgetful functor $F:\Mm^{\Cc} \to \Mm_A$ is separable, hence $\rho_M$ also has a left inverse $f_1 \in \Mm^{\Cc}$ (see \cite{CMZ}). So $f_1$ is split epimorphism in $\Mm^{\Cc}$.\\
The epimorphism
$${\rm can}:\ \Sigma^* \ot_T \Sigma\to \Cc$$
is split in $\Mm_A$, hence it is also split in $\Mm^\Cc$, because $\Cc$ is coseparable.
Tut then $f_2=M\ot_A{\rm can}$ is a split epimorphism in $\Mm^\Cc$, hence
$g=f_1\circ f_2$ is a split epimorphism in $\Mm^\Cc$.\\
Now $P = \Ker(g) \in \Mm^{\Cc}$. Taking $M = P$ in the above reasoning, we obtain another split epimorphism in $\Mm^\Cc$:
$$h:\ P \ot_A (\Sigma^* \ot_T \Sigma) \to P$$
From $g$ split epi and the natural transformation $\zeta$, we obtain the following commutative diagram with exact rows
$$\begin{diagram}
0 & \rTo & P & \rTo & M\ot_A (\Sigma^*\ot_T\Sigma) & \rTo^g & M & \rTo & 0 \\
&&\uTo_{\zeta_P} & & \uTo_{\zeta_{M\ot_A (\Sigma^*\ot_T\Sigma)}} && \uTo_{\zeta_M} && \\
0 & \rTo & FG(P) & \rTo & FG(M\ot_A (\Sigma^*\ot_T\Sigma)) & \rTo^{FG(g)} & FG(M) & \rTo & 0
\end{diagram}$$
We have a similar diagram for $h$:
$$\begin{diagram}
0 & \rTo & \Ker(h) & \rTo & P\ot_A (\Sigma^*\ot_T\Sigma) & \rTo^h & P & \rTo & 0 \\
&&\uTo_{\zeta_{Ker(h)}} & & \uTo_{\zeta_{P\ot_A (\Sigma^*\ot_T\Sigma)}} && \uTo_{\zeta_P} && \\
0 & \rTo & FG(\Ker(h)) & \rTo & FG(P\ot_A (\Sigma^*\ot_T\Sigma)) & \rTo^{FG(h)} & FG(P) & \rTo & 0
\end{diagram}$$
With these two commutative diagrams with exact rows, we can make a third one
$$\begin{diagram}
P\ot_A \Sigma^* \ot_T \Sigma & \rTo & M\ot_A (\Sigma^*\ot_T\Sigma) & \rTo^g & M & \rTo & 0 \\
\uTo_{\zeta_{P\ot_A \Sigma^* \ot_T \Sigma}} & & \uTo_{\zeta_{M\ot_A (\Sigma^*\ot_T\Sigma)}} && \uTo_{\zeta_M} && \\
FG(P\ot_A \Sigma^* \ot_T \Sigma) & \rTo & FG(M\ot_A (\Sigma^*\ot_T\Sigma)) & \rTo^g & FG(M) & \rTo & 0
\end{diagram}$$
Ty \leref{5.3} the first two vertical arrows are isomorphisms. 
From the lemma of 5, it now follows that $\zeta_M$ is an isomorphism.
\end{proof}

We have an inverse to \prref{5.5}. But first, let us give a characterization of the
coseparability of the comatrix coring.

\begin{lemma}\lelabel{5.6}
Let $\Sigma\in {}_T\Mm$ and $\Sigma^*\in \Mm_T$ be totally faithful. Then the 
comatrix coring $\Dd = \Sigma^*\ot_T \Sigma$ is coseparable if and only if
the map $l: \ \Sigma\ot^{\Dd}\Sigma^* \to \Sigma\ot_A\Sigma^*$ is split monomorphism in ${}_T\Mm_T$. 
\end{lemma}

\begin{proof}
Let $\Dd$ be coseparable. By \leref{5.4}, there exists a $\Sigma$-normalized
$\theta\in V_2$. By \leref{5.t} we have a surjective projection $t:\ \Sigma \ot_A \Sigma^* \to \Sigma \ot^{\Dd} \Sigma^*$ in ${}_T\Mm_T$ and so we have that $t\circ l = I_{\Sigma \ot^{\Dd} \Sigma^*}$. \\
Conversely, if $l$ is split mono in ${}_T\Mm_T$, then there exists a $t:\ \Sigma \ot_A \Sigma^* \to \Sigma \ot^{\Dd} \Sigma^*$ in ${}_T\Mm_T$ such that $t\circ l = I_{\Sigma \ot^{\Dd} \Sigma^*}$. Now denote by $i$ the composition of $t$ and the canonical injection of $\Sigma \ot^{\Dd}\Sigma^*$ into $T = \Sigma \ot^{\Cc}\Sigma^*$ and define $\theta = \varepsilon_{\Dd} \circ (I_{\Sigma^*} \ot_T i \ot_T I_{\Sigma})$. 
It is clear that $\theta$ is $(A\hbox{-}A)$-bilinear. Let us describe $\theta$ more explicitely
$$
\theta((g \ot_B u)\ot_A (h \ot_B v)) = \varepsilon_{\Dd}(g.t(u\ot_A h) \ot_B v)  g(t(u\ot_A h)v).
$$
On one hand we have
\begin{eqnarray*}
&&\hspace*{-2cm}
\sum_i g\ot_B e_i\theta((f_i\ot_B u)\ot_A (h \ot_B v))\\
 &=& \sum_i g\ot_B e_i(f_i(t(u \ot_A h)v)) =
 g\ot_B t(u \ot_A h)v,
\end{eqnarray*}
while on the other hand 
\begin{eqnarray*}
&&\hspace*{-2cm}
\sum_i \theta(g\ot_B u\ot_A h \ot_B e_i)f_i \ot_B v\\ &=& \sum_i g(t(u\ot_A h)e_i)f_i \ot_B v =
 g\ot_B t(u \ot_A h)v
\end{eqnarray*}
This proves that $\theta \in \tilde{V}_2$ ($\tilde{V}_2$ is defined as $V_2$ but with $\Cc$ replaced by $\Dd$). Finally
\begin{eqnarray*}
&&\hspace*{-2cm}
\theta(\Delta_{Dd}(u\ot_B g)) = \theta(u \ot_B e \ot_B g) =
 \varepsilon_{\Dd}(u.t(e) \ot_B g) \\
&=& \varepsilon_{\Dd}(u.e \ot_B g) =
 \varepsilon_{\Dd}(u \ot_B g)
\end{eqnarray*}
This concludes the proof that $\Dd$ is coseparable.
\end{proof}

Combining the results of \prref{5.5}, \leref{5.6} and \prref{3.3}, we obtain

\begin{theorem}\thlabel{5.6}
Let $\Cc$ be projective as a right $A$-module, $T = \End^{\Cc}(\Sigma) \cong \Sigma\ot^{\Cc}\Sigma^*$ and
$\Dd = \Sigma^*\ot_B\Sigma$. The following are equivalent
\begin{enumerate}
\item
\begin{itemize}
\item $(F,G)$ is a pair of inverse equivalences;
\item $l: \Sigma\ot^{\Dd}\Sigma^* \to \Sigma\ot_A\Sigma^*$ is split mono in ${}_T\Mm_T$.
\end{itemize}
\item
\begin{itemize}
\item $\Cc$ is a coseparable coring;
\item ${\rm can}$ is surjective.
\end{itemize}
\end{enumerate}
\end{theorem}

\thref{5.6} is also new in the situation where $\Sigma=A$, with right $\Cc$-coaction
$\rho(a)=xa$ with $x\in G(\Cc)$ a grouplike element. Then \thref{5.6} takes the following
form.

\begin{corollary}\colabel{5.7}
Let $(\Cc,x)$ be an $A$-coring with a fixed grouplike element, and
$$T=A^{{\rm co}\Cc}=\{b\in A~|~xb=bx\}.$$
If $\Cc$ is projective as a right $A$-module, then the following assertions are equivalent.
\begin{enumerate}
\item
\begin{itemize}
\item $(\bullet\ot_T A,(\bullet)^{{\rm co}\Cc})$ is a pair of inverse equivalences;
\item $i:\ T\to A$ is split mono in ${}_T\Mm_T$.
\end{itemize}
\item
\begin{itemize}
\item $\Cc$ is a coseparable coring;
\item ${\rm can}:\ A\ot_T A\to \Cc$, ${\rm can}(a\ot_Tb)=axb$, is surjective.
\end{itemize}
\end{enumerate}
\end{corollary}

\section{Frobenius corings}\selabel{6}
Recall that an $A$-coring $\Cc$ is called Frobenius if the right adjoint of the forgetful functor $\Mm^\Cc\to\Mm_A$ is also a left adjoint. The forgetful functor and its
adjoint are then called a Frobenius pair. $\Cc$ is Frobenius if and only if
$\Cc\in {}_A\Mm$ is locally projective and
there exists a  bijective map $j\in{}_A\Hom_{\*C}(\Cc,\*C)$. In this situation,
$\Cc$ is finitely generated and projective as a left and right $A$-module,
and the categories ${_A\Mm^\Cc}$ and ${_A\Mm_{\*C}}$ are isomorphic. 
$\Cc$ is Frobenius if and only if there exists a Frobenius system, consisting of
a pair $(z,\theta)$, with
$z\in \Cc^A=\{c\in\Cc~|~ac=ca, {\rm ~for~all~}a\in A\}$ and $\theta\in{_A\Hom_A}(\Cc\otimes_A\Cc,A)$ such that the following conditions hold:
\begin{itemize}
\item $c_{(1)}\theta(c_{(2)}\otimes_Ad)=\theta(c\otimes_Ad_{(1)})d_{(2)}$, for all $c, d \in \Cc$,
\item $\theta(z\otimes_A c)=\theta(c\otimes_A z)=\varepsilon_\Cc(c)$, for all $c\in \Cc$.
\end{itemize}
For details we refer to \cite{CMZ}.\\
One implication of
our next result is a generalization of \cite[Theorem 2.7]{CVW}.

\begin{proposition}\prlabel{6.1}
Let $\Cc$ be an $A$-coring, take $\Sigma\in\Mm^\Cc_{fgp}$ and consider the Morita context $\CC$ associated to $\Sigma$ as introduced in \seref{4}. If $\Cc$ is Frobenius, then there exists an isomorphism of $(A,B)$-bimodules $J:\ \Sigma^*\to Q$. The Morita context $\CC$ is isomorphic to the Morita context $\widetilde{\CC}=(T,\*C,\Sigma,\Sigma^*,\mu,\tau)$, where the left $\*C$-action on $\Sigma^*$ and the maps $\mu$ and $\tau$ are given explicitly by
$$(g\cdot f)(u)=\theta\Bigl(z_{(1)}g(z_{(2)})\otimes_Af(u_{[0]})u_{[1]}\Bigr),$$
$$\mu:\ \Sigma^*\otimes_T\Sigma\to \*C,~~\mu(f\otimes_T u)(c)=\theta(c\otimes_Af(u_{[0]})u_{[1]}),$$
$$\tau:\ \Sigma\otimes_{\*C}\Sigma^*\to T,~~
\tau(u\otimes_{\*C} f)(v)=u_{[0]}\theta(u_{[1]}\otimes_Af(v_{[0]})v_{[1]}),$$
where $f\in\Sigma^*,~g\in\*C,~u,v\in\Sigma,~c\in\Cc$ and ($z$,$\theta$) is a Frobenius system for $\Cc$.

Conversely, if $\Cc$ satisfies the equivalent conditions of \thref{4.11}, and if $\Sigma^*$ and
$Q$ are isomorphic as $(A,B)$-bimodules, then $\Cc$ is a Frobenius coring.
\end{proposition}

\begin{proof}
Applying the functor $G=\Hom^\Cc(\Sigma,\bullet)$ to the Frobenius map
$j\in{_A\Hom_{\*C}}(\Cc,\*C)$, we obtain the following isomorphism in $\Mm_B$
$$\Hom^\Cc(\Sigma,j):\Hom^\Cc(\Sigma,\Cc)\to\Hom^\Cc(\Sigma,\*C).$$
Now $\Hom^\Cc(\Sigma,\Cc)\cong\Hom_A(\Sigma,A)=\Sigma^*$ and $\Hom_{\*C}(\Sigma,\*C)\cong{^\Cc\Hom}(\Cc,\Sigma^*)=Q$ (see \prref{4.8}). Hence we obtain
an isomorphism of right $B$-modules $J:\ \Sigma^*\to Q=\Hom^\Cc(\Sigma,\*C)$. A straightforward
computation shows that $J$ is given by the formula
$$J(f)(u)=j(f(u_{[0]})u_{[1]}),$$
for all $f\in \Sigma^*$ and $u\in\Sigma$. Let us show that $J$ is left $A$-linear.
$$
J(af)(u)=j((af)(u_{[0]})u_{[1]})=j(a(f(u_{[0]})u_{[1]}))=
aj(f(u_{[0]})u_{[1]})=aJ(f)(u),$$
where we just used the $A$-linearity of $\Sigma^*$ and $j$. The left
$\*C$-module structure on $Q$ to $\Sigma^*$:
\begin{eqnarray*}
&&\hspace*{-2cm}
(g\cdot f)(u)=J^{-1}(g\cdot J(f))(u)=\varepsilon\circ j^{-1}((g\cdot J(f))(u))\\
&=&\varepsilon\circ j^{-1}(g\# J(f)(u))= \varepsilon\circ j^{-1}(g\# j(f(u_{[0]})u_{[1]}))\\
&=&\varepsilon(z\cdot(g\# j(f(u_{[0]})u_{[1]})))=
\varepsilon(z_{(1)}\cdot(j(f(u_{[0]})u_{[1]})(z_{(2)}g(z_{(2)}))))\\
&=&\varepsilon(z_{(1)}\theta(z_{(2)}g(z_{(3)})\otimes_Af(u_{[0]})u_{[1]}))\\
&=&\varepsilon(z_{(1)})\theta(z_{(2)}g(z_{(3)})\otimes_Af(u_{[0]})u_{[1]})\\
&=&\theta(z_{(1)}g(z_{(2)})\otimes_Af(u_{[0]})u_{[1]}).
\end{eqnarray*}
Now take the Morita context $\CC$ from \equref{4.5.1}. Using the isomorphisms $J$
and $\beta:\ \Hom_{\*C}(\Sigma,\*C)\to{^\Cc\Hom}(\Cc,\Sigma^*)=Q$ from \prref{4.8},
we find the connecting maps of the Morita context: 
\begin{eqnarray*}
&&\hspace*{-2cm}
\mu(f\otimes_T u)(c)=\beta(J(f))(c)(u)=J(f)(u)(c)\\
&=&j(f(u_{[0]})u_{[1]})(c)=\theta(c\otimes_Af(u_{[0]})u_{[1]}).
\end{eqnarray*}
Remark that $\mu(f\otimes_T u)=J(f)(u)$.
\begin{eqnarray*}
&&\hspace*{-2cm}
\tau(u\otimes_{\*C} f)(v)=u_{[0]}(\beta(J(f))(u_{[1]})(v))=u_{[0]}(J(f)(v)(u_{[1]}))\\
&=&u_{[0]}j(f(v_{[0]})v_{[1]})(u_{[1]})=u_{[0]}\theta(u_{[1]}\otimes_Af(v_{[0]})v_{[1]}).
\end{eqnarray*}
Conversely, assume that $\Cc$ satisfies the equivalent conditions of \thref{4.11}, and let
$J:\ \Sigma^*\to Q$ be an isomorphism of  $(A,B)$-bimodules. 
Then we find that $J\ot_B\Sigma: \Sigma^*\ot_B\Sigma$,
$\mu:\ Q\ot_B\Sigma\to {}^*\Cc$ and ${\rm can}:\ \Sigma^*\ot_B\Sigma \to \Cc$
are $(A,{}^*\Cc)$-bimodule isomorphisms. Then $j=\mu\circ(J\ot_B\Sigma)\circ
{\rm can}^{-1}$ is an $(A,{}^*\Cc)$-bimodule isomorphism between $\Cc$ and ${}^*\Cc$,
and we find that ${}^*\Cc$ is Frobenius.
\end{proof}

Recall from \cite{Br2} that a coring $\Cc$ is called right coFrobenius if $\Cc\in{}_A\Mm$
is locally projective, and if there exists a injective $j\in{_A\Hom_{\*C}}(\Cc,\*C)$.
Observe that this notion is not left-right symmetric. The following result may be viewed
as a generalization of \cite[Theorem 2.10]{BDR}.

\begin{proposition}\prlabel{6.2}
Let $\Cc$ be an $A$-coring, take $\Sigma\in\Mm^\Cc_{\rm fgp}$ and consider the Morita context $\CC$ associated to $\Sigma$ as in \seref{4}. If $\Cc$ is right coFrobenius, then there exists a monomorphism of $(A,B)$-bimodules $J:\ \Sigma^*\to Q$.\\
Conversely, if the equivalent conditions of \thref{4.15} are satisfied, and if there
exists a monomorphism of $(A,B)$-bimodules $J:\ \Sigma^*\to Q$, then
$\Cc$ is right coFrobenius.
\end{proposition}

\begin{proof}
We construct $J$ in the same way as in \prref{6.1}, namely $J(f)(u)=
j(f(u_{[0]})u_{[1]})$. Let us show that $J$ is injective. If $J(f(u))=0$, then
$j(f(u_{[0]})u_{[1]})=0$, for all $u\in \Sigma$. Since $j$ is injective, this implies that
$f(u_{[0]})u_{[1]}=0$, hence $0=\varepsilon(f(u_{[0]})u_{[1]})=f(u_{[0]})\varepsilon(u_{[1]})
=f(u_{[0]}\varepsilon(u_{[1]}))=f(u)$, for all $u\in \Sigma$, and $f=0$.
Conversely, if $J:\ \Sigma^*\to Q$ is a monomorphism of $(A,B)$-modules, then we have
a morphism of $(A,{}^*\Cc)$-modules
$J\ot_B\Sigma:\ \Sigma^*\ot_B\Sigma\to Q\ot_B\Sigma$, which is injective since
$\Sigma\in {}_B\Mm$ is flat. \thref{4.15} also tells us that ${\rm can}:\ \Sigma^*\ot_B\Sigma\cong\Cc$
and $\mu:\ Q\ot_B\Sigma\to ({}^*\Cc)^{\rm rat}$ are isomorphisms. Then
$j=\mu\circ(J\ot_B\Sigma)\circ
{\rm can}^{-1}:\ \Cc\to ({}^*\Cc)^{\rm rat}\subset {}^*\Cc$ is a monomorphism of
$(A,{}^*\Cc)$-bimodules.
\end{proof}

Before we state our next results, we need some Lemmas.

\begin{lemma}
Let $\Cc$ be an $A$-coring. Then
\lelabel{6.3}
$${_A\Hom_{\*C}(\Cc,\*C)}\cong {_{\Cc^*}\Hom_A(\Cc,\Cc^*)}$$
\end{lemma}

\begin{proof}
Take and $(A,{}^*\Cc)$-bimodule map $j:\ \Cc\to {}^*\Cc$. Then
$\tilde{j}:\ \Cc\to\Cc^*$ defined by $\tilde{j}(c)(d)=j(d)(c)$ is a $(\Cc^*,A)$-bimodule map.
All verifications are straightforward.
\end{proof}

For $P\in {}_A\Mm$, and $R\subset P$, we define
$$R^\bot=\{p\in P~|~r(p)=0,\forall r\in R\}.$$
Recall that a ring $A$ is called Pseudo-Frobenius ring (or PF ring) if
$A$ is an injective cogenerator of $\Mm_A$. Examples of such rings are symmetric algebras, Frobenius algebras, quasi-Frobenius rings  or QF rings and finite dimensional Hopf algebras. Moreover, if $R$ is a principal ideal domain, then $R/I$ is a PF ring (and even a QF ring) for every ideal $I$. If $R$ is a QF-ring, then ${\rm M}_n(R)$ is also a QF ring.

\begin{lemma}\lelabel{6.4}
Let $A$ be a PF ring, $\Cc$ an $A$-coring, and $j:\ \Cc\to {}^*\Cc$
an $(A,{}^*\Cc)$-bimodule map.
 With notation as in \leref{6.3},
 $\tilde{\j}$ is an injection if and only if $\im j$ is dense in the finite topology on $\*C$. In this case $(\*C)^{\rm rat}$ is dense.
\end{lemma}

\begin{proof}
The injectivity of $\tilde{j}$ is equivalent to
$$\tilde{j}(d)(c)=j(d)(c)=0,~~{\rm for~all}~d\in \Cc~~\Longrightarrow~~c=0$$
and to
$$f(c)=0,~~{\rm for~all}~f\in \im(j)~~\Longrightarrow~~c=0,$$
which can be restated as follows:
$$\im j^\perp=\{0\}$$
Since $A$ is a PF ring, this is equivalent to $\im j$ is dense in the finite topology on $\*C$,
see \cite[Theorem 1.8]{A2}. Finally, $\im j\subset (\*C)^{\rm rat}$.
\end{proof}

\begin{lemma}\lelabel{6.5}
Let $A$ be a commutative ring and $\Cc$ an $A$-coalgebra. If $(\Cc^*)^{\rm rat}$ is dense in the finite topology on $\Cc^*$, then ${_{\Cc^*}\Hom}((\Cc^*)^{\rm rat},M)={_{\Cc^*}\Hom}(\Cc^*,M)$ for every $M\in {^\Cc_{\rm fgp}\Mm}$.
\end{lemma}

\begin{proof}
Let $E(M)$ be the injective envelope of $M\in {}^\Cc\Mm$. Then by 
\cite[9.5]{BrzezinskiWisbauer} $E(M)$ is also injective as a left $\Cc^*$-module and we can extend any left $\Cc^*$-linear $\chi:\ (\Cc^*)^{\rm rat}\to M\subset E(M)$ to $\bar{\chi}:\Cc^*\to E(M)$.

Since $(\Cc^*)^{\rm rat}$ is dense, it has left local units on $M$, so we can take $e\in(\Cc^*)^{\rm rat}$ such that $e\cdot m=m$ with $m=\bar{\chi}(\varepsilon_\Cc)$. We find $e\cdot \bar{\chi}(\varepsilon_\Cc) =\bar{\chi}(e\#\varepsilon_\Cc)=\bar{\chi}(e)\in M$. Furthermore, for any $f\in\Cc^*$, $\bar{\chi}(f)=\bar{\chi}(f\#\varepsilon_\Cc)=f\cdot \bar{\chi}(\varepsilon_\Cc)\in M$, so $\bar{\chi}\in{_{\Cc^*}\Hom}(\Cc^*,M)$.

Finally, $\bar{\chi}$ is unique: suppose that there exists a $\xi\in{_{\Cc^*}\Hom}(\Cc^*,M)$ which has also the property that $\xi(f)=\chi(f)$ for all $f\in(\Cc^*)^{\rm rat}$, and take a local unit $e\in(\Cc^*)^{\rm rat}$ for $(\xi-\bar{\chi})(\varepsilon_\Cc)$. Then
$$(\xi-\bar{\chi})(\varepsilon_\Cc)
=e\cdot(\xi-\bar{\chi})(\varepsilon_\Cc)=(\xi-\bar{\chi})(e)=0.$$
Consequently
$$(\xi-\bar{\chi})(f)=(\xi-\bar{\chi})(f\#\varepsilon_\Cc)=
f\cdot(\xi-\bar{\chi})(\varepsilon_\Cc)=0,$$
finishing the proof.
\end{proof}

\begin{proposition}\prlabel{6.6}
Let $A$ be a commutative ring, and
$\Cc$ an $A$-coalgebra. Take $\Sigma\in\Mm^\Cc_{\rm fgp}$ and consider the Morita context $\CC$ associated to $\Sigma$ as in \seref{4}. If $A$ is a commutative PF ring and $\Cc$ is left and right coFrobenius, then there exists a monomorphism of $(A,B)$-bimodules $J:\ \Sigma^*\to Q$ and an epimorphism of $(A,B)$-bimodules $J':\ \Sigma^*\to Q$. 
Conversely, if $A$ is a commutative PF ring, and
if the equivalent conditions of \thref{4.15} are satisfied, then the existence
of $J$ and $J'$ as implies that $\Cc$ is left and right coFrobenius.
\end{proposition}

\begin{proof}
The monomorphism $J$ is constructed as in \prref{6.2}.\\
Since $A$ is a $PF$ ring and therefore injective in $\Mm_A$, $\Cc$ is injective in $\Mm^\Cc$. By \cite[9.5]{BrzezinskiWisbauer}, $\Cc$ is also injective as a $\Cc^*$-module,  so the injective left coFrobenius morphism $j':\ \Cc\to\Cc^{*\rm rat}$ splits, and $\Cc$ is a direct summand of $\Cc^{*\rm rat}$ as a $\Cc^*$-module. We obtain an epimorphism
\begin{equation}\eqlabel{6.6.1}
{_{\Cc^*}\Hom}(\Cc^{*\rm rat},\Sigma^*)\to {_{\Cc^*}\Hom}(\Cc,\Sigma^*).
\end{equation}
$\Cc$ is right coFrobenius, so it follows from \leref{6.4} that
 $(\*C)^{\rm rat}$ is dense in the finite topology. From \leref{6.5}, it follows that
 ${_{\Cc^*}\Hom}(\Cc^{*\rm rat},\Sigma^*)=\Sigma^*$.\\
To prove the converse, we proceed as in \prref{6.2}. The existence of
the monomorphism $J$ implies that $\Cc$ is right coFrobenius.
Using the fact that ${\rm can}$ and $\mu$ are isomorphisms, we find a
$(B,A)$-bimodule epimorphism
$$\tilde{j}'=\mu\circ (J'\ot_B\Sigma)\circ{\rm can}^{-1}:\
\Cc\to\Sigma^*\otimes_B\Sigma\to Q\otimes_B\Sigma\to(\*C)^{\rm rat}$$
Since $(\*C)^{\rm rat}$ is dense, the dual morphism is defined on $\Cc$,
and is injective by \leref{6.4}. So  we find that $\Cc$ is also left
coFrobenius.
\end{proof}

\section{The case where $\Sigma=\Cc$}\selabel{7}
Let $\Cc$ be an $A$-coring which is finitely generated and projective as
a right $A$-module. $\Cc$ is a right $\Cc$-comodule, and,
by \prref{1.1}, $\Cc^*$ is a left $\Cc$-comodule. Consider the pairs of adjoint
functors $(F,G)$ and $(F',G')$
 from \prref{1.4}, where we take $\Sigma=\Cc$ and $B=T=\End^\Cc(\Cc)\cong\Cc^*$:
 \begin{equation}\eqlabel{7.1.1}
F:\Mm_{\Cc^*}\to\Mm^\Cc,~~ F(N)=N\otimes_{\Cc^*}\Cc,
\end{equation}
\begin{equation}\eqlabel{7.1.2}
G:\Mm^\Cc\to\Mm_{\Cc^*},~~ G(M)=\Hom^\Cc(\Cc,M)\cong M\otimes^\Cc\Cc^*,
\end{equation}
and
$$F':{_{\Cc^*}\Mm}\to{^\Cc\Mm},~~ F'(N)=\Cc^*\otimes_{\Cc^*}N,$$
$$G':{^\Cc\Mm}\to{_{\Cc^*}\Mm},~~ G'(M)={^\Cc\Hom}(\Cc^*,M)\cong\Cc\otimes^\Cc M.$$
Since $\Cc$ is finitely generated and projective as a right $A$-module,
the categories ${_{\Cc^*}\Mm}$ and ${^\Cc\Mm}$ are isomorphic. The isomorphism
and its inverse are given by the functors $F'$ and $G'$.\\
The associated comatrix coring is $\Dd=\Cc^*\otimes_{\Cc^*}\Cc$ and the cannonical map 
$${\rm can}:\Dd\to \Cc,{\rm can}(f\otimes_{\Cc^*} c)=f(c_{(1)})c_{(2)}$$
is the canonical isomorphism. We also have two Morita contexts. The first context
is the one from \reref{4.5}  (2), with $M=\Cc$. We find
$$\CC'=(T=\Cc^*,\Cc^*,Q=\End^\Cc(\Cc)\cong\Cc^*,\Cc^*,\tau,\mu),$$
with $\tau=\mu$ the canonical isomorphism $\Cc^*\ot_{\Cc^*}\Cc^*\to \Cc^*$.
This Morita context is the trivial one connecting $\Cc^*$ to itself.\\
The second context is the one from \reref{4.5}  (3), with $\Sigma=\Cc$.
This leads us to
\begin{equation}\eqlabel{7.1.3}
\CC=(T=\Cc^*,{}^*\Cc,\Cc,Q={}^\Cc\Hom(\Cc,\Cc^*),\tau,\mu).
\end{equation}
We now want to investigate when $(F,G)$ is a pair of inverse equivalences.
In the situation where $\Cc$ is also finitely generated and projective as a
left $A$-module, the answer is given by \thref{4.11}. We obtain the following
result. 

\begin{corollary}\colabel{7.1}
Let $\Cc$ be an $A$-coring which is finitely generated and projective as
a left and right $A$-module. Then the following assertions are equivalent.
\begin{enumerate}
\item $\Cc\in {}_{\Cc^*}\Mm$ is faithfully flat;
\item $\Cc\in {}_{\Cc^*}\Mm$ is a progenerator;
\item the Morita context \equref{7.1.3} is strict;
\item $(F,G)$ from (\ref{eq:7.1.1}-\ref{eq:7.1.2}) is a pair of inverse
equivalences.
\end{enumerate}
\end{corollary}

We will now give other sufficient conditions for $(F,G)$ to be a pair of inverse
equivalences. Recall first that $M\in \Mm^\Cc$ is called right $\Cc$-coflat
if it is flat as a right $A$-module, and if $M\ot^{\Cc}-:\ {}^\Cc\Mm\to \dul{\rm Ab}$
is exact. A similar definition applies to left $\Cc$-comodules.

\begin{lemma}
\lelabel{7.2} Let $A$ be a ring.
With $M \in {\Mm}^{\Cc}, N \in {}^{\Cc}{\Mm}_A$ and $P \in {}_A{\Mm}$, the natural map
$$f:\ (M\ot^{\Cc} N)\ot_A P\to M\ot^{\Cc} (N\ot_A P)$$
is an isomorphism in each of the following situations:
\begin{enumerate}
\item $P\in {}_A\Mm$ is flat;
\item $M\in \Mm^\Cc$ is coflat.
\end{enumerate}
\end{lemma}

\begin{proof}
1. Recall that $M\ot^{\Cc} N$ is defined by the exact sequence
$$0\to{}M\ot^{\Cc} N\to{}M\ot_A N\rightrightarrows M\ot_A {\Cc}\ot_A N$$
Using the fact that $P$ is $A$-flat, we obtain a commutative diagram
with exact rows
$$\begin{matrix}
0&\to{}&(M\ot^{\Cc} N)\ot_A P&\to{}&M\ot_A N\ot_A P&\rightrightarrows&
M\ot_A {\Cc}\ot_A N\ot_A P\cr
&&\downarrow{f}&&\downarrow{\cong}&&\downarrow{\cong}\cr
0&\to{}&M\ot^{\Cc} (N\ot P)&\to{}&M\ot_A N\ot P&\rightrightarrows&
M\ot_A {\Cc}\ot_A N\ot P\cr\end{matrix}$$
and the result follows from the Five Lemma.\\

2. Recall the definition of the tensor product over $\ZZ$:
$N\ot_\ZZ P=\ZZ(N\times P)/I$, where $I$ is the ideal generated by elements
of the form
$$(n,p+q)-(n,p)-(n,q)~~;~~(n+m,p)-(n,p)-(m,p)~~;~~
(nx,p)-(n,xp)$$
This means we can construct an exact sequence of left ${\Cc}$-comodules
$$0\to{}J'\to{}\ZZ(N\times P)/{I'}\to{}N\ot P\to{}0$$
where $I'$ is the ideal generated by elements of the form $(n+m,p)-(n,p)-(m,p)$, and $J'$ the ideal in $\ZZ(N\times P)/{I'}$ that is generated by elements of the form 
$$\ol{(n,p+q)-(n,p)-(n,q)}~~;~~\ol{(nx,p)-(n,xp)}$$
Now, using the right ${\Cc}$-coflatness of $M$, we find
a commutative diagram with exact rows
$$\begin{matrix}
0&\to{}&M\ot^{\Cc} J'&\to{}&M\ot^{\Cc}\ZZ(N\times P)/{I'}&\to{}&
M\ot^{\Cc}(N\ot_\ZZ P)&\to{}&0\cr
&&\downarrow{=}&&\downarrow{\cong}&&\uparrow{f}&&\cr
0&\to{}&J''&\to{}&\ZZ((M\ot^{\Cc} N)\times P)/{I''}&\to{}&
(M\ot^{\Cc} N)\ot_\ZZ P&\to{}&0\cr\end{matrix}$$
and it follows the from Five Lemma that 
$(M\ot^{\Cc} N)\ot_A P\cong M\ot^{\Cc} (N\ot_A P)$. We then obtain the following
commutative diagram with exact rows
$$\begin{matrix}
M\ot^{\Cc} (N\ot_\ZZ A\ot_\ZZ  P) & \rightrightarrows & M\ot^{\Cc} (N\ot_\ZZ P) 
& \rightarrow & M\ot^{\Cc} (N \ot_A P) & \rightarrow & 0 \cr
\downarrow{\cong}&&\downarrow{\cong}&&\downarrow{}&& \cr
(M\ot^{\Cc} N)\ot_\ZZ A\ot_\ZZ  P & \rightrightarrows & (M\ot^{\Cc} N)\ot_\ZZ P 
& \rightarrow & (M\ot^{\Cc} N) \ot_A P & \rightarrow & 0 \cr
\cr
\end{matrix}$$
The second row is the defining exact sequence of the tensor product
$(M\ot^\Cc N)\ot_A P$. The result then follows from the Lemma of 5.
\end{proof}

\begin{theorem}\thlabel{7.3}
Let $\Cc$ be an $A$-coring which is finitely generated and projective as
a right $A$-module. If $\Cc\in {_{\Cc^*}\Mm}$ is flat
and $\Cc^*\in {^\Cc\Mm}$ is coflat, then the adjoint pair $(F,G)$
from (\ref{eq:7.1.1}-\ref{eq:7.1.2}) is a pair of inverse equivalences.
\end{theorem}

\begin{proof}
We first prove that the counit of the adjunction is an isomorphism. The counit is 
given by the formula
$$\zeta_M:\ (M\otimes^\Cc\Cc^*)\otimes_{\Cc^*}\Cc\to M,~~
\zeta((\sum_im_i\otimes_Af_i)\otimes_{\Cc^*}c)=\sum_im_if_i(c)$$
By \leref{7.2}, we have isomorphisms
$$(M\otimes^\Cc\Cc^*)\otimes_{\Cc^*}\Cc\cong M\otimes^\Cc(\Cc^*\otimes_{\Cc^*}\Cc)
\cong M\otimes^\Cc\Cc\cong M$$
The composition of these isomorphisms is precisely$\zeta_M$, hence
$\zeta_M$ is an isomorphism.\\
Similar arguments show that the unit of the adjunction is an isomorphism.
Let $\sum_i e_i\ot_A f_i\in \Cc\ot_A\Cc^*$ be a dual basis for $\Cc\in\Mm_A$.
Then the unit is given by the formula 
$$\nu_N:\ N\to(N\otimes_{\Cc^*}\Cc)\otimes^\Cc\Cc^*,~~\nu_N(n)=(n\otimes_{\Cc^*}e_i)\otimes_Af_i.$$
Observe that 
\begin{equation}\eqlabel{7.3.1}
\rho_{\Cc^*}(\varepsilon_\Cc)=\varepsilon_\Cc(e_{i(1)})e_{i(2)}\otimes_Af_i=e_i\otimes_Af_i.
\end{equation}
If $\Cc^*$ is coflat as left $\Cc$-comodule, then we have the following
$$N\cong N\otimes_{\Cc^*}\Cc^* \cong N\otimes_{\Cc^*}(\Cc\otimes^\Cc\Cc^*) \cong (N\otimes_{\Cc^*}\Cc)\otimes^\Cc\Cc^*,$$
Using \equref{7.3.1}, we find that this composition is $\nu_N$, hence
$\nu_N$ is an isomorphism.
\end{proof}

\begin{theorem}\thlabel{7.4}
Let $\Cc$ be a Frobenius $A$-coring. Then the adjoint pair $(F,G)$
from (\ref{eq:7.1.1}-\ref{eq:7.1.2}) is a pair of inverse equivalences.
\end{theorem}

\begin{proof}
Recall that a Frobenius coring is finitely generated projective on both sides.
By \coref{7.1}, it suffices to show that the Morita context $\CC$ is strict.\\
Take a Frobenius isomorphism $j:\ \Cc\to\Cc^*$ in ${_{\Cc^*}\Mm_A}$.
We first prove that the map $\tau$ of the Morita context $\CC$ is surjective. The map is given explicitly by 
$$\tau:\ \Cc\otimes_{\*C}Q\to\Cc^*,~~\tau(c\otimes q)=q(c).$$ 
Observe that $j\in Q$, and consider the map 
$$\tau':\ \Cc^*\to\Cc\otimes_{\*C}Q,~~\tau'(f)=j^{-1}(f)\otimes j$$
$\tau'$ is a right inverse of $\tau$, so $\tau$ is surjective and a fortiori bijective.\\
It follows from \prref{4.8} that
$$Q={^\Cc\Hom(\Cc,\Cc^*)}\cong \Hom^\Cc(\Cc,\*C)=\tilde{Q}.$$
Now $\mu$ is given by
$$\mu:\ \tilde{Q}\otimes_{\Cc^*}\Cc\to\*C,~~\mu(\tilde{q}\otimes c)=\tilde{q}(c).$$ 
the inverse of $\mu$ is the map $\mu':\ \*C\to\tilde{Q}\otimes_{\Cc^*}\Cc,~~\mu'(f)=\tilde{j}\otimes\tilde{j}^{-1}(f)$.
\end{proof}

We consider again the case $\Sigma=\Cc\in\Mm^\Cc_{\rm fgp}$, but this time we take $B=A$
instead of $B=T$. The map $\ell:\ B=A\to T=\End^\Cc(\Cc)=\Cc^*$, is now the usual ring homomorphism $i:\ A\to \Cc^*$, given by $i(a)(c)=a\varepsilon_\Cc(c)$. 
We have the two following pairs of adjoint functors $(F,G)$ and $(F',G')$.
$$F:\ \Mm_A\to \Mm^\Cc,~~F(N)=N\ot_A\Cc$$
$$G:\ \Mm^\Cc\to\Mm_A,~~G(M)=M\ot^\Cc \Cc^*$$
$$F':\ {}_A\Mm\to {}^\Cc\Mm,~~F'(N)=\Cc^*\ot_AN$$
$$G':\ {}^\Cc\Mm\to {}_A\Mm,~~G'(M)=\Cc\ot^\Cc M\cong M$$
$G'$ is the forgetful functor. Now we know that the functor $F$ also has
a left adjoint $H$, and that the forgetful functor $G'$ has a right
adjoint $H'=\Cc\ot_A\bullet$. Now recall that the coring $\Cc$ is called a
Frobenius coring if the forgetful functor $\Mm^\Cc\to \Mm^A$ is Frobenius,
which means that it has a right adjoint which is at the same time a left
adjoint. This is equivalent to the forgetful functor $G'$ being Frobenius,
see \cite[Theorem 35]{CMZ}; more equivalent conditions are given in
\cite[Theorem 36]{CMZ}. Using the adjoint pairs $(F,G)$ and $(F',G')$,
we can state more equivalent conditions:

\begin{proposition}\prlabel{7.5}
Let $\Cc$ be an $A$-coring which is finitely generated and projective as
a right $A$-module. With notation as above, the following assertions are
equivalent:
\begin{enumerate}
\item $\Cc$ is a Frobenius coring;
\item $G$ is isomorphic to the forgetful functor $H$;
\item $F'$ is (isomorphic to) the functor $H'=\Cc\otimes_A\bullet$; 
\item $G$ is a left adjoint of $F$;
\item $G'$ is a left adjoint of $F'$.
\end{enumerate}
\end{proposition} 
 
\begin{center}
{\bf Acknowledgment}
\end{center}
We thank Bachuki Mesablishvili
for helping us with \seref{2}.

\end{document}